\documentclass[12pt,a4paper]{article}
\usepackage{amstext}
\usepackage{fancyhdr}
\usepackage{amsfonts,graphicx,bezier, amssymb}
\usepackage{amsmath}
\usepackage{caption}
\usepackage{mathtools}
\usepackage{mathrsfs}
\usepackage[utf8]{inputenc}
\usepackage{tikz}
\usepackage{rotating}
\usepackage[multiple]{footmisc}
\newenvironment{prf}{\noindent{\bf{Proof:}}~~}{\hfill\rule{1ex}{1ex}\vskip1.5ex}
\newcommand{\Z}{\mathbb Z}
\newcommand{\Q}{\mathbb Q}

\newcommand{\beqa}{\begin{eqnarray}}
	\newcommand{\enqa}{\end{eqnarray}}
\newcommand{\beq}{\begin{eqnarray*}}
	\newcommand{\enq}{\end{eqnarray*}}

\newtheorem{cor}{Corollary}[section]
\newtheorem{prop}{Proposition}[section]
\newtheorem{defn}{Definition}[section]
\newtheorem{exam}{{\bf Example}}[section]
\newtheorem{thm}{Theorem}[section]

\newtheorem{lem}{Lemma}[section]

\newcommand{\noi}{\noindent}

\makeatletter
\providecommand*{\twoheadrightarrowfill@}{%
	\arrowfill@\relbar\relbar\twoheadrightarrow
}
\providecommand*{\twoheadleftarrowfill@}{%
	\arrowfill@\twoheadleftarrow\relbar\relbar
}
\providecommand*{\xtwoheadrightarrow}[2][]{%
	\ext@arrow 0579\twoheadrightarrowfill@{#1}{#2}%
}
\providecommand*{\xtwoheadleftarrow}[2][]{%
	\ext@arrow 5097\twoheadleftarrowfill@{#1}{#2}%
}
\makeatother
\usepackage[para]{manyfoot}

\begin{document}
 \begin{center}
	{\bf \Large Modules (co)reduced relative to another module}
\end{center}
\begin{center}
	Tilahun Abebaw\footnote{Department of Mathematics, Addis Ababa University, P.O. Box 1176, Addis Ababa, Ethiopia \newline email: tilahun.abebaw@aau.edu.et, amanuel.mamo@aau.edu.et, zelalem.teshome@aau.edu.et}, Amanuel Mamo\footnotemark[1]{}\footnote{This work forms part of the second author's PhD thesis.}, 
	David Ssevviiri\footnote {Department of Mathematics, Makerere University, P.O. Box 7062, Kampala Uganda, email: david.ssevviiri@mak.ac.ug}\footnote{Corresponding author} and
	Zelalem Teshome\footnotemark[1]	
\end{center}

\begin{abstract}\noi  Let $R$ be a commutative unital ring, $\mathfrak{ a}$ an ideal of $R$ and $M$ a fixed $R$-module. We introduce and study generalisations of $\mathfrak{a}$-reduced modules, $\mathfrak{R}_{\mathfrak{ a}}$  and $\mathfrak{a}$-coreduced modules, $\mathfrak{C}_{\mathfrak{ a}}$ studied in \cite{kyomuhangi2020locally, Annet-David : generalised reduced, David-Application-I, ssevviiri: App II}. They are called modules $\mathfrak{ a}$-reduced with respect to $M$, $\mathfrak{R}^{M}_{\mathfrak{ a}}$ and modules $\mathfrak{ a}$-coreduced with respect to $M$, $\mathfrak{C}^{M}_{\mathfrak{ a}}$.  The paper lifts several known results about $\mathfrak{R}_{\mathfrak{ a}}$ and $\mathfrak{C}_{\mathfrak{ a}}$ in \cite{kyomuhangi2020locally, Annet-David : generalised reduced, David-Application-I, ssevviiri: App II} to larger classes of modules $\mathfrak{R}^{M}_{\mathfrak{ a}}$ and $\mathfrak{C}^{M}_{\mathfrak{ a}}$ respectively. Results about $\mathfrak{R}_{\mathfrak{ a}}$ and $\mathfrak{C}_{\mathfrak{ a}}$ can be retrieved by taking $M= R$.
\end{abstract}
\vspace*{0.4cm}\\
{\bf MSC 2020} Mathematics Subject Classification: 13C13,13D07, 13D45, 18A35. 
\section{Introduction}
\begin{paragraph}\noi
Let $R$ be a commutative unital ring, and let $\mathfrak{ a}$ be an ideal of  $R$. The functors
 $\Gamma_{\mathfrak{ a}}(-):=\underset{k}\varinjlim\text{Hom}_{R}(R/\mathfrak{ a}^{k},-)$ and $\Lambda_{ \mathfrak{a}}(-):= \underset{k}\varprojlim R/\mathfrak{ a}^{k}\otimes_{R}-$ on the category of $R$-modules, $R\text{-Mod}$ which are respectively called the \textit{$\mathfrak{ a}$-torsion functor} and the \textit{$\mathfrak{ a}$-adic completion functor} are dual to each other. They have been used to study different phenomena both in commutative algebra and algebraic geometry \cite{Alonso-Lipman,brodmann2013local,K.Divaani, Hartshorne-Local cohomology,Hartshorne-conjecture, Annet's dissertation,kyomuhangi2020locally, Annet-David : generalised reduced,David-Application-I, ssevviiri: App II}.  An $R$-module $M$ is \textit{$\mathfrak{ a}$-reduced} (resp. \textit{$\mathfrak{ a}$-coreduced})  if $\mathfrak{ a}^{2}m=0$ implies $\mathfrak{ a}m=0$  for all $m\in M$ (resp. $\mathfrak{ a}^{2}M=\mathfrak{ a}M$). It was shown in \cite[Proposition 2.2]{David-Application-I} (resp. \cite[Proposition 2.3]{David-Application-I}) that $M$ is $\mathfrak{ a}$-reduced (resp. $\mathfrak{ a}$-coreduced) if and only if $\mathfrak{a}\Gamma_{\mathfrak{ a}}(M)=0$ (resp. $\mathfrak{ a}\Lambda_{ \mathfrak{a}}(M)=0$). The modules $\mathfrak{a}\Gamma_{\mathfrak{ a}}(M)$ (resp. $\mathfrak{ a}\Lambda_{ \mathfrak{a}}(M)$) are therefore obstructions to $M$ being $\mathfrak{ a}$-reduced (resp. $\mathfrak{ a}$-coreduced). A module is \textit{reduced} (resp. \textit{coreduced}) if it is $\mathfrak{ a}$-reduced (resp. $\mathfrak{ a}$-coreduced) for all ideals $\mathfrak{ a}$ of $R$. The objective of the papers \cite{kyomuhangi2020locally, David-Application-I, ssevviiri: App II} was in part to make the study of reduced modules and coreduced modules which were introduced in \cite{Lee and Zhou- reduced modules} and \cite{Ansari-Farshadadifar} respectively categorical. This offers some advantages over studying them using element-wise methods. For instance, it provides avenues to prove possible equivalences between them and other module categories, and to also know about their homological properties. This paper achieves something more; being able to generalise notions to larger classes of modules.
 \end{paragraph}
 
\begin{paragraph}\noi
 Given a fixed $R$-module $M$, we define new classes of $R$-modules called \textit{$\mathfrak{ a}$-reduced modules with respect to $M$}, and \textit{$\mathfrak{ a}$-coreduced modules with respect to $M$} which strictly contain $\mathfrak{ a}$-reduced modules,   $\mathfrak{R}_{\mathfrak{ a}}$ and $\mathfrak{ a}$-coreduced modules, $\mathfrak{ C}_{\mathfrak{ a}}$ respectively, see Proposition \ref{Prop: reduced w.r.t M} (resp. Proposition \ref{sec1:Prop- strongly I-coreduced}). 	
 We say that an $R$-module $N$ is \textit{$\mathfrak{ a}$-reduced with respect to $M$} (resp. \textit{ $\mathfrak{ a}$-coreduced with respect to $M$}) if an $R$-module $\text{Hom}_{R}(M,N)$ is $\mathfrak{ a}$-reduced (resp. $M\otimes_{R}N$ is $\mathfrak{ a}$-coreduced). It turns out that $N$ is $\mathfrak{ a}$-reduced with respect to $M$ (resp. $\mathfrak{ a}$-coreduced with respect to $M$) if and only if $\mathfrak{ a}\Gamma_{ \mathfrak{a}}(M,N)=0$ (resp. $\mathfrak{ a}\Lambda_{ \mathfrak{a}}(M,N)=0$), see Proposition \ref{Prop: reduced w.r.t M} (resp. Proposition \ref{sec1:Prop- strongly I-coreduced}),  	 where $\Gamma_{ \mathfrak{a}}(M,-):=\underset{k}{\varinjlim}\text{Hom}_{R}(M/\mathfrak{ a}^{k}M,-)$ and $\Lambda_{ \mathfrak{a}}(M,-):=\underset{k}\varprojlim(M/\mathfrak{ a}^{k}M\otimes_{R}-)$ are the well-known generalised $\mathfrak{ a}$-torsion functor and the generalised $\mathfrak{ a}$-adic completion functor respectively; see for instance \cite{T T Nam-GLH for artinian modules, T.T. Nam-generalised I-adic completion} where they were utilised.
  This paper lifts several results proved in \cite{kyomuhangi2020locally, David-Application-I, ssevviiri: App II} about $\mathfrak{R}_{\mathfrak{ a}}$ and $\mathfrak{ C}_{\mathfrak{ a}}$ to larger classes of modules $\mathfrak{R}_{\mathfrak{ a}}^{M}$ and $\mathfrak{ C}_{\mathfrak{ a}}^{M}$ respectively. The results obtained in \cite{kyomuhangi2020locally,Annet-David : generalised reduced, David-Application-I, ssevviiri: App II} can be retrieved by taking $M=R$. 
\end{paragraph}

\begin{paragraph}\noi The paper has seven sections. In the first section, we have the introduction. In section two (resp. section three)  the classes $\mathfrak{R}_{\mathfrak{ a}}^{M}$ (resp. $\mathfrak{ C}_{\mathfrak{ a}}^{M}$) have been introduced and their basic properties given. These properties are utilised in the subsequent sections. In section four, we show that $\mathfrak{R}_{\mathfrak{ a}}^{M}$ (resp. $\mathfrak{ C}_{\mathfrak{ a}}^{M}$) is a complete subcategory of $R\text{-Mod}$ and forms a pretorsion-free class of modules (resp. a cocomplete subcategory of $R\text{-Mod}$ and forms a pretorsion class of modules). In section five, we give duality theorems. For instance,  it shown that the functor $\Gamma_{ \mathfrak{a}}(M,-): \mathfrak{R}^{M}_{\mathfrak{ a}}\rightarrow \mathfrak{ C}^{M}_{\mathfrak{ a}}$ is right adjoint to the functor $\Lambda_{ \mathfrak{a}}(M,-):  \mathfrak{ C}^{M}_{\mathfrak{ a}}\rightarrow \mathfrak{R}^{M}_{\mathfrak{ a}}$  generalising the Greenlees-May Duality given in \cite[Theorem 3.4]{David-Application-I}. Section six contains applications of $\mathfrak{R}^{M}_{\mathfrak{ a}}$ and $\mathfrak{ C}^{M}_{\mathfrak{ a}}$ to the generalised local cohomology modules $\text{H}^{i}_{\mathfrak{a}}(M,N)$ and the generalised local homology modules  $\text{H}_{i}^{\mathfrak{ a}}(M,N)$ respectively. For instance, if $R$ is Noetherian and $N\in\mathfrak{R}^{M}_{\mathfrak{ a}}$ (resp. $N\in\mathfrak{C}^{M}_{\mathfrak{ a}}$), then so is $\text{H}^{i}_{\mathfrak{a}}(M,N)$ (resp. $\text{H}^{i}_{\mathfrak{a}}(M,N)$). The same is true if Noetherian is replaced with Artinian. In the last section, section seven, further applications of $\mathfrak{R}^{M}_{\mathfrak{ a}}$ and $\mathfrak{C}^{M}_{\mathfrak{ a}}$ to the generalised local cohomology modules and to the generalised local homology modules via spectral sequences are given. They determine when $\text{H}^{i}_{\mathfrak{a}}(M,N)\in \mathfrak{R}^{M}_{\mathfrak{ a}}$ (resp. $\text{H}_{i}^{\mathfrak{a}}(M,N)\in \mathfrak{C}^{M}_{\mathfrak{ a}}$). Furthermore, they respectively facilitate the vanishing of  $\text{H}^{p}_{\mathfrak{ a}}(M,\text{H}^{q}_{\mathfrak{ a}}(M,N))$ and $\text{H}^{\mathfrak{ a}}_{p}(M,\text{H}^{\mathfrak{ a}}_{q}(M,N))$ whenever either $p\ne0$ or $q\ne0$.
	\end{paragraph}

\begin{paragraph}\noi
\textbf{Notation}. Throughout this paper, $R$ is a commutative unital ring and $\mathfrak{ a}$ is an ideal of $R$. $R\text{-Mod}$ denotes the category of all $R$-modules and $R\text{-mod}$ is the full subcategory of $R\text{-Mod}$ consisting of finitely generated $R$-modules. Let $M\in R\text{-Mod}$. The full subcategories of $R\text{-Mod}$ which consist of $\mathfrak{ a}$-reduced $R$-modules, $\mathfrak{ a}$-coreduced $R$-modules,  $R$-modules $\mathfrak{ a}$-reduced with respect to $M$, and $R$-modules $\mathfrak{ a}$-coreduced with respect to $M$ are denoted by $\mathfrak{R}_{\mathfrak{ a}}, \mathfrak{C}_{\mathfrak{ a}}, \mathfrak{R}_{\mathfrak{ a}}^{M}$ and $\mathfrak{C}_{\mathfrak{ a}}^{M}$ respectively.
\end{paragraph}

\section{Reduced with respect to a module}

\begin{paragraph}\noi
Let $\mathfrak{ a}$ be an ideal of a ring $R$, and let $M\in R\text{-Mod}$.
In this section, we introduce and study the basic properties of $R$-modules which are $\mathfrak{a}$-reduced  with respect to $M$.
\end{paragraph}

\begin{defn}\label{Defn: reduced w.r.t M}
{\normalfont	An $R$-module $N$ is}  $\mathfrak{ a}$-reduced with respect to $M$ {\normalfont if $\text{Hom}_{R}(M,N)$ is an $\mathfrak{ a}$-reduced $R$-module.} {\normalfont  We call $N$} reduced with respect to $M$ {\normalfont if it is $\mathfrak{ a}$-reduced with respect to $M$ for all ideals $\mathfrak{ a}$ of $R$.} 
\end{defn}

\begin{paragraph}\noi
An $R$-module is $\mathfrak{ a}$-reduced if and only if it is $\mathfrak{ a}$-reduced with respect to  $R$.
\end{paragraph}

\begin{lem}\label{Isomorphisms: The Gammas and the Lamdas}
For all $R$-modules $M$ and $N$, there are isomorphisms of $R$-modules {\normalfont  $\Gamma_{\mathfrak{ a}}(M,N)\cong  \Gamma_{\mathfrak{ a}}(\text{Hom}_{R}(M,N))$} and {\normalfont $\Lambda_{\mathfrak{ a}}(M,N)\cong \Lambda_{\mathfrak{ a}}(M\otimes_{R}N).$ }
\end{lem}

\begin{prf}
By definition, $\Gamma_{\mathfrak{ a}}(M,N):=\underset{k}{\varinjlim}\text{Hom}_{R}(M/\mathfrak{ a}^{k}M,N)\cong \underset{k}{\varinjlim}\text{Hom}_{R}(R/\mathfrak{ a}^{k}\otimes_{R}M, N) $$\cong\underset{k}{\varinjlim}\text{Hom}_{R}(R/\mathfrak{ a}^{k}, \text{Hom}_{R}(M,N))\cong \Gamma_{\mathfrak{ a}}(\text{Hom}_{R}(M,N)).$ Dually,
$ \Lambda_{\mathfrak{ a}}(M,N):=\underset{k}{\varprojlim}(M/\mathfrak{ a}^{k}M\otimes_{R}N)\cong \underset{k}{\varprojlim}(R/\mathfrak{ a}^{k}\otimes_{R}M)\otimes_{R} N $$\cong\underset{k}{\varprojlim}(R/\mathfrak{ a}^{k}\otimes_{R}(M\otimes_{R}N))\cong \Lambda_{\mathfrak{ a}}(M\otimes_{R}N).$	
\end{prf}

\begin{paragraph}\noi Lemma \ref{Isomorphisms: The Gammas and the Lamdas} demonstrates that the functors $\Gamma_{ \mathfrak{a}}(M,-)$ (resp. $\Lambda_{ \mathfrak{a}}(M,-)$) are nothing but compositions of $\text{Hom}_{R}(M,-)$ and $\Gamma_{ \mathfrak{a}}(-)$ (resp.  $M\otimes_{R}- $ and $\Lambda_{ \mathfrak{a}}(-)$).
\end{paragraph}

\begin{lem} \label{isomorphism: gener. gamma. version}
If $M$ and $N$ are $R$-modules with $M$ finitely presented, then there is an isomorphism {\normalfont $\Gamma_{ \mathfrak{a}}(M,N)\cong \text{Hom}_{R}(M,\Gamma_{ \mathfrak{a}}(N))$}. 
\end{lem}

\begin{prf}
$\Gamma_{\mathfrak{ a}}(M,N):=\underset{k}{\varinjlim}\text{Hom}_{R}(M/\mathfrak{ a}^{k}M,N)\cong \underset{k}{\varinjlim}\text{Hom}_{R}(M\otimes_{R}R/\mathfrak{ a}^{k}, N)$. This together with the Hom-Tensor adjunction and the fact that Hom commutes with the direct limits in the second variable since $M$ is finitely presented, we have
$\Gamma_{\mathfrak{ a}}(M,N)\cong\underset{k}{\varinjlim}\text{Hom}_{R}(M, \text{Hom}_{R}(R/\mathfrak{ a}^{k},N))\cong \text{Hom}_{R}(M,\underset{k}{\varinjlim}\text{Hom}_{R}(R/\mathfrak{ a}^{k},N))$$\cong \text{Hom}_{R}(M,\Gamma_{ \mathfrak{a}}(N)).$ 
\end{prf}

\begin{paragraph}\noi From Lemma \ref{Isomorphisms: The Gammas and the Lamdas} and Lemma \ref {isomorphism: gener. gamma. version}, we conclude that if $M$ is finitely presented, then the functors $\Gamma_{ \mathfrak{a}}(-)$ and $\text{Hom}_{R}(M, -)$ commute.
\end{paragraph}

\begin{prop}\label{Prop: Gammas for red w.r.t M}
Let {\normalfont $M,N\in R\text{-Mod}$} with $M$ finitely presented. Then $N\in\mathfrak{R}^{M}_{\mathfrak{ a}}$ if and only if $\Gamma_{ \mathfrak{a}}(N) \in\mathfrak{R}^{M}_{\mathfrak{ a}}$.
\end{prop}

\begin{prf}
In light of \cite[Corollary 2.13]{David-Application-I} and Lemma \ref{isomorphism: gener. gamma. version}, an $R$-module $\text{Hom}_{R}(M,N)$ is $\mathfrak{ a}$-reduced  if and only if $\Gamma_{ \mathfrak{a}}(\text{Hom}_{R}(M,N))\cong \Gamma_{ \mathfrak{a}}(M,N)\cong \text{Hom}_{R}(M,\Gamma_{ \mathfrak{a}}(N))$ is $\mathfrak{ a}$-reduced.
\end{prf}
\begin{prop}\label{Prop: reduced w.r.t M}
Let $\mathfrak{ a}$ be an ideal of a ring $R$. The following statements are equivalent for every pair of $R$-modules $M$ and $N$.
{\normalfont \begin{itemize}
\item[$(1)$] $N$ is $\mathfrak{ a}$-reduced with respect to $M$,
\item[$(2)$] $\text{Hom}_{R}(M/\mathfrak{a}M, N) \cong \text{Hom}_{R}(M/\mathfrak{a}^{2}M, N),$

\item[$(3)$] $\Gamma_\mathfrak{a}(M,N)\cong  \text{Hom}_{R}(M/\mathfrak{a}M, N),$ 
\item[$(4)$] $\mathfrak{ a}\Gamma_{\mathfrak{ a}}(M,N)= 0,$
\item[$(5)$] $\Gamma_{ \mathfrak{a}}(M,N)$ is $\mathfrak{ a}$-reduced.
\end{itemize}}
\end{prop}
\begin{prf}
\begin{itemize}
\item[$(1)\Rightarrow (2)$] Since $N\in\mathfrak{R}^{M}_{\mathfrak{ a}}$, by Definition \ref{Defn: reduced w.r.t M}, $ \text{Hom}_{R}(M,N)\in\mathfrak{R}^{}_{\mathfrak{ a}}$. So, by \cite [Proposition 2.2]{David-Application-I}, $\text{Hom}_{R}(R/\mathfrak{ a},\text{Hom}_{R}(M,N))\cong \text{Hom}_{R}(R/\mathfrak{ a}^{2},\text{Hom}_{R}(M,N)).$  
It follows, by the Hom-Tensor adjunction and the isomorphism $M/\mathfrak{ a}^{k}M\cong R/\mathfrak{ a}^{k}\otimes_{R} M$ for all $ k\in \Z^{+}$, that
$ \text{Hom}_{R}(M/\mathfrak{a}M, N) \cong \text{Hom}_{R}(M/\mathfrak{a}^{2}M, N).$
\item[$(2)\Rightarrow (3)$] From the hypothesis, $ \text{Hom}_{R}(M/\mathfrak{a}^{k}M, N)\cong  \text{Hom}_{R}(M/\mathfrak{a}M, N)~ \text{for all} ~ k\in\Z^{+}.$ Thus, $\Gamma_{\mathfrak{a}}(M,N):= \underset{k}\varinjlim~\text{Hom}_{R}(M/\mathfrak{a}^{k}M, N)\cong \text{Hom}_{R}(M/\mathfrak{a}M, N)$.
\item[$(3)\Rightarrow (4)$] This follows, since $\Gamma_\mathfrak{a}(M,N)\cong \text{Hom}_{R}(M/\mathfrak{a}M,N)$ implies that $\mathfrak{ a}\Gamma_\mathfrak{a}(M,N)\cong  \text{Hom}_{R}(\mathfrak{ a}M/\mathfrak{ a}M,N))=0$.
\item[$(4)\Rightarrow (5)$] By Lemma \ref{Isomorphisms: The Gammas and the Lamdas}, $\Gamma_\mathfrak{a}(M,N)\cong \Gamma_\mathfrak{a}(\text{Hom}_{R}(M, N)).$ This implies that $\mathfrak{ a}\Gamma_{ \mathfrak{a}}(M,N)\cong \mathfrak{ a}\Gamma_{ \mathfrak{a}}(\text{Hom}_{R}(M,N)).$ Since $\mathfrak{ a}\Gamma_{ \mathfrak{a}}(M,N)=0$ by hypothesis, $\mathfrak{ a}\Gamma_{ \mathfrak{a}}(\text{Hom}_{R}(M,N))=0$. Thus, $\text{Hom}_{R}(M,N)$ is $\mathfrak{ a}$-reduced \cite[Proposition 2.2]{David-Application-I}, and  this is true if and only if $\Gamma_{ \mathfrak{a}}(\text{Hom}_{R}(M,N))\cong \Gamma_{\mathfrak{ a}}(M,N)$ is $\mathfrak{ a}$-reduced \cite[Corollary 2.13]{David-Application-I}.
\item[$(5)\Rightarrow (1)$] By \cite[Corollary 2.13]{David-Application-I} and Lemma \ref{Isomorphisms: The Gammas and the Lamdas}, $\Gamma_{ \mathfrak{a}}(M,N)\cong \Gamma_{ \mathfrak{a}}(\text{Hom}_{R}(M,N))$ is $\mathfrak{ a}$-reduced if and only if $\text{Hom}_{R}(M,N)$ is $\mathfrak{ a}$-reduced. This proves that $N\in\mathfrak{R}^{M}_{\mathfrak{ a}}$.
\end{itemize}
\end{prf}
\begin{lem}\label{Lem: further example-red.w.r.t M}
Let $M$ and $ N$ be $R$-modules.  
\begin{itemize}
\item [$(1)$] If $N\in\mathfrak{R}^{}_{\mathfrak{ a}}$, then $N\in\mathfrak{R}^{K}_{\mathfrak{ a}}$ for any {\normalfont $K\in R\text{-Mod}$}, i.e.,  {\normalfont $\mathfrak{R}_{\mathfrak{ a}}\subseteq \mathfrak{R}^{K}_{\mathfrak{ a}}$ for any $K\in R\text{-Mod}$}.
\item[$(2)$] If $M\in\mathfrak{ C}_{\mathfrak{ a}}$, then {\normalfont $\mathfrak{R}^{M}_{\mathfrak{ a}}=R\text{-Mod}$}.
\end{itemize} 	 
\end{lem}
\begin{prf}
\begin{itemize}
\item[$(1)$] Let $N\in\mathfrak{R}_{\mathfrak{ a}}$. We show that $\text{Hom}_{R}(K,N)\in \mathfrak{R}_{\mathfrak{ a}}$ for any $K\in R\text{-Mod}$.  Suppose $\varphi \in \text{Hom}_{R}(K,N)$ and $a\in \mathfrak{ a}$ such that $a^{2}\varphi=0.$ So, $a^{2}\varphi(m)=0$ for all $m\in K$. Since $\varphi(m)\in N$ and $N\in \mathfrak{R}_{\mathfrak{ a}}$ by hypothesis, $a\varphi(m)=0$ for all $m$, and thus $a\varphi=0$. This shows that $\text{Hom}_{R}(K,N)\in\mathfrak{R}_{\mathfrak{ a}}$. Therefore, $N\in\mathfrak{R}^{K}_{\mathfrak{ a}}$. 
\item[$(2)$]$\mathfrak{ a}M=\mathfrak{ a}^{2}M$, by hypothesis. So,
$\text{Hom}_{R}(M/\mathfrak{ a}M,K)\cong\text{ Hom}_{R}(M/\mathfrak{ a}^{2}M, K)$ for every $
K\in R\text{-Mod}$. The conclusion is clear by Proposition \ref{Prop: reduced w.r.t M}. 
\end{itemize} 
\end{prf}
\begin{paragraph}\noi 
Lemma \ref{Lem: further example-red.w.r.t M} shows that $\mathfrak{R}_{\mathfrak{ a}}$ is contained in {\normalfont $\mathfrak{R}_{\mathfrak{a}}^{M}$}.
However, not all $\mathfrak{a}$-reduced modules with respect to a given {\normalfont $M\in R\text{-Mod}$} are $\mathfrak{ a}$-reduced, as the following examples illustrate. 
\begin{enumerate} 
\item[$(1)$]	Let $\Z$ be the ring of integers, and $\Z_{n}$ the group of integers modulo $n$. Since $\Z_{2}$ is $(2)$-coreduced as a $\Z$-module, the $\Z$-module $ \text{Hom}_{\Z}(\Z_{2}, \Z_{4})$ is $(2)$-reduced \cite[Proposition 2.6 (1)]{David-Application-I}. So, by  Definition \ref{Defn: reduced w.r.t M}, $\Z_{4}\in\mathfrak{R}^{\Z_{2}}_{(2)}$. But, $\Z_{4}\notin\mathfrak{R}_{(2)}$, since there exists $2\in (2)$ and $\bar{3}\in\Z_{4}$ for which  $2.\bar{3}=\bar{2}\ne0$ but $2^{2}\bar{3}=0$.
\item[$(2)$] The $\Z$-module $\Q/\Z$ is not reduced \cite[page 5]{David-Application-I}. However, $\Q/\Z$ is reduced with respect any coreduced module $M$, see Lemma \ref{Lem: further example-red.w.r.t M} (2). In particular, $\Q/\Z$ is reduced with respect to itself.
	\end{enumerate} 
\end{paragraph}
\section{Coreduced with respect to a module}
\begin{paragraph}\noi	Let $\mathfrak{ a}$ be an ideal of a ring $R$, and let $M\in R\text{-Mod}$.
	In this section, we introduce and study the basic properties of $R$-modules which are $\mathfrak{a}$-coreduced with respect to $M$. 
	\end{paragraph}
\begin{defn}\label{Defn: Cored w.r.t M}
{\normalfont	An $R$-module $N$ is} $\mathfrak{ a}$-coreduced with respect to $M$ {\normalfont if an $R$-module $M\otimes_{R}N$ is $\mathfrak{ a}$-coreduced.} {\normalfont If $N$ is $\mathfrak{ a}$-coreduced with respect to $M$ for all ideals $\mathfrak{ a}$ of $R$, then $N$ is called} coreduced with respect to $M$.
\end{defn}
\begin{paragraph}\noi
An $R$-module is $\mathfrak{ a}$-coreduced if and only if it is $\mathfrak{ a}$-coreduced with respect to  $R$.
\end{paragraph}
\begin{prop}\label{sec1:Prop- strongly I-coreduced}
Let $\mathfrak{ a}$ be an ideal of  $R$, and let $M,N$ be $R$-modules. The following statements are equivalent:
\end{prop}
\begin{itemize}
\item[$(1)$] $N$ is $\mathfrak{ a}$-coreduced with respect to $M$, 
\item[$(2)$] $M/\mathfrak{a}M\otimes_{R} N \cong M/\mathfrak{a}^{2}M \otimes_{R} N,$
\item[$(3)$] $\Lambda_{\mathfrak{a}}(M,N)\cong M/\mathfrak{a}M\otimes_{R} N,$ 
\item[$(4)$]$\mathfrak{a}\Lambda_{\mathfrak{a}}(M,N)= 0 $,
\item[$(5)$] $\Lambda_{ \mathfrak{a}}(M,N)$ is $\mathfrak{ a}$-coreduced.
\end{itemize}
\begin{prf}
\begin{itemize}
\item[$(1)\Rightarrow (2)$] Suppose that $N\in\mathfrak{ C}^{M}_{\mathfrak{ a}}$. By Definition \ref{Defn: Cored w.r.t M}, $ M\otimes_{R}N$ is $\mathfrak{a}$-coreduced. So, $
R/\mathfrak{ a}\otimes_{R} (M\otimes_{R} N)\cong R/\mathfrak{ a}^{2}\otimes_{R} (M\otimes_{R} N)$ \cite [Proposition 2.3]{David-Application-I}. Applying associativity of the tensor product, we get $
(R/\mathfrak{ a}\otimes_{R} M)\otimes_{R} N\cong (R/\mathfrak{ a}^{2}\otimes_{R} M)\otimes_{R} N.$ Therefore,  $M/\mathfrak{a}M\otimes_{R} N \cong M/\mathfrak{a}^{2}M \otimes_{R} N.$   
\item[$(2)\Rightarrow (3)$] From $(2)$, it is readily seen that $M/\mathfrak{ a}M\otimes_{R} N\cong M/\mathfrak{ a}^{k}M\otimes_{R} N$ for all $k\in \mathbb{Z}^{+}$. Thus, $\Lambda_{\mathfrak{a}}(M,N):= \underset{k}\varprojlim~(M/\mathfrak{a}^{k}M\otimes_{R} N)\cong \underset{k}\varprojlim~(M/\mathfrak{a}M\otimes_{R} N)\cong M/\mathfrak{a}M\otimes_{R}N$.
\item[$(3)\Rightarrow (4)$] Since $\Lambda_{ \mathfrak{a}} (M,N)\cong M/\mathfrak{ a}M\otimes_{R}N$, the conclusion $\mathfrak{ a}\Lambda_{ \mathfrak{a}}(M, N)=0$ is apparent. 
\item[$(4)\Rightarrow (5)$] Since, by hypothesis, $\mathfrak{ a}\Lambda_{ \mathfrak{a}} (M,N)=0$ and by  Lemma \ref{Isomorphisms: The Gammas and the Lamdas},  $\Lambda_{ \mathfrak{a}} (M,N)\cong \Lambda_{ \mathfrak{a}}(M\otimes_{R}N)$, $\mathfrak{ a}\Lambda_{ \mathfrak{a}}(M\otimes_{R}N)\cong 0$. So, $M\otimes_{R}N$ is $\mathfrak{ a}$-coreduced  \cite[Proposition 2.3]{David-Application-I}. This is true if and only if $\mathfrak{a}^{2}(M\otimes_{R}N)=\mathfrak{a}(M\otimes_{R}N)$ if and only if $\Lambda_{ \mathfrak{a}}(\mathfrak{ a}^{2}(M\otimes_{R}N))=\Lambda_{ \mathfrak{a}}(\mathfrak{ a}(M\otimes_{R}N))$ if and only if $\mathfrak{ a}^{2} \Lambda_{ \mathfrak{a}}(M\otimes_{R}N)=\mathfrak{ a} \Lambda_{ \mathfrak{a}}(M\otimes_{R}N)$ if and only if $\mathfrak{ a}^{2} \Lambda_{ \mathfrak{a}}(M,N)=\mathfrak{ a} \Lambda_{ \mathfrak{a}}(M,N)$. This proves the claim.
\item[$(5)\Rightarrow (1)$] If $\Lambda_{ \mathfrak{a}}(M,N)$ is $\mathfrak{ a}$-coreduced, then $\mathfrak{ a}^{2} \Lambda_{ \mathfrak{a}}(M,N)=\mathfrak{ a} \Lambda_{ \mathfrak{a}}(M,N)$. This holds true if and only if $M\otimes_{R}N$ is $\mathfrak{ a}$-coreduced, see proof of $(4)\Rightarrow (5)$. Thus, $N\in\mathfrak{C}^{M}_{\mathfrak{ a}}$.
\end{itemize}
\end{prf}
\begin{exam}
 If  the ideal $\mathfrak{ a}$ is idempotent, then {\normalfont $\mathfrak{R}_{\mathfrak{ a}}= \mathfrak{R}_{\mathfrak{ a}}^{M}=\mathfrak{C}_{\mathfrak{ a}}=\mathfrak{C}_{\mathfrak{ a}}^{M}=R\text{-Mod}$}.
\end{exam}
\begin{prop}\label{Lem: cor tensor with cor is cor} Suppose that $M$ and $N$ are $R$-modules. If either $M$ or $N$ is $\mathfrak{ a}$-coreduced, then so is the $R$-module $M\otimes_{R}N$. In particular, {\normalfont $\mathfrak{C}_{\mathfrak{ a}}\subseteq \mathfrak{C}_{\mathfrak{ a}}^{M}$ for any $M$, and $\mathfrak{ C}^{M}_{\mathfrak{ a}}=R\text{-Mod}$ whenever $M$ is $\mathfrak{ a}$-coreduced}. 
\end{prop}
\begin{prf}
Suppose that either $aM= a^{2}M$ or $a^{2}N=aN$ for all $a\in\mathfrak{ a}$. Then either $aM\otimes_{R}N\cong a^{2}M\otimes_{R}N ~\text{or}~ M\otimes_{R}aN\cong M\otimes_{R}a^{2}N~\text{for all}~ a\in\mathfrak{ a}.$ In each case,  $a(M\otimes_{R}N)\cong a^{2}(M\otimes_{R}N)$ for every $a\in \mathfrak{ a}$. Therefore, $M\otimes_{R}N$ is $\mathfrak{ a}$-coreduced.
\end{prf}
\begin{paragraph}\noi Proposition \ref{Lem: cor tensor with cor is cor} says that $\mathfrak{C}_{\mathfrak{ a}}$ is contained in $\mathfrak{C}_{\mathfrak{ a}}^{M}$. The following examples demonstrate that the containment is strict.
	\begin{enumerate}  
		\item[$(1)$] Since $\Z_{2}$ over $\Z$ is $(2)$-coreduced, $ \Z_{2}\otimes_{\Z} \Z_{4}$ is $(2)$-coreduced [Proposition \ref{Lem: cor tensor with cor is cor}]. By  Definition \ref{Defn: Cored w.r.t M}, $\Z_{4}$ is $(2)$-coreduced with respect to $\Z_{2}$. However, $\Z_{4}$ is not $(2)$-coreduced; $(2)^{2}~\Z_{4}=0$ but $(2)~ \Z_{4}=\{\bar{0},\bar{2}\}\ne 0$.
		\item[$(2)$] The $\Z$-module $\Z$ is not coreduced. 
		However, $\Z$ is coreduced with respect to any coreduced module $M$, see Proposition \ref{Lem: cor tensor with cor is cor}. In particular, $\Z$ is coreduced with respect to $\Q$.
		\end{enumerate}
\end{paragraph}

\begin{prop}\label{Prop: Lemma 2.2 of our first paper} Let $M$ and $N$ be $R$-modules. If {\normalfont $N \in \mathfrak{R}^{M}_{\mathfrak{ a}}$} (resp. {\normalfont $N\in\mathfrak{C}^{M}_{\mathfrak{ a}}$}), then  {\normalfont $\text{Hom}_{R}(M,N)\in\mathfrak{R}^{M}_{\mathfrak{ a}}$} (resp. {\normalfont $M\otimes_{R}N\in \mathfrak{C}^{M}_{\mathfrak{ a}}$}). 
\end{prop}
\begin{prf}
Let $M$ be an $R$-module, and let $N\in\mathfrak{R}^{M}_{\mathfrak{ a}}$. By Definition \ref{Defn: reduced w.r.t M}, $\text{Hom}_{R}(M,N)$ is $\mathfrak{ a}$-reduced. By Lemma \ref{Lem: further example-red.w.r.t M}, $\text{Hom}_{R}(M,N)\in\mathfrak{R}^{M}_{\mathfrak{ a}}$ follows. This asserts the first statement.   
 If $N\in\mathfrak{C}^{M}_{\mathfrak{ a}}$, then $M\otimes_{R}N\in\mathfrak{ C}_{\mathfrak{ a}}$, and thus $M\otimes_{R}N\in\mathfrak{C}^{M}_{\mathfrak{ a}}$ [Definition \ref{Defn: Cored w.r.t M} and Proposition \ref{Lem: cor tensor with cor is cor}].
\end{prf}
   
\begin{prop}\label{Hom and tensor is str. (co)reduced}
Let $M$ be an $R$-module. If $X$ and $Y$ are $R$-modules with {\normalfont $X\in\mathfrak{C}^{M}_{\mathfrak{ a}}$}, then 
  {\normalfont $\text{Hom}_{R}(X,Y)\in\mathfrak{R}^{M}_{\mathfrak{ a}}$} and {\normalfont $X\otimes_{R}Y\in \mathfrak{C}^{M}_{\mathfrak{ a}}$}. 
\end{prop}
\begin{prf}
By Proposition \ref{sec1:Prop- strongly I-coreduced}, $M/\mathfrak{ a}M\otimes_{R}X\cong M/\mathfrak{ a}^{2}M\otimes_{R}X$ for all $X\in\mathfrak{C}^{M}_{\mathfrak{a}}$ . Applying the Hom-Tensor adjunction,
{\normalfont 
$\text{Hom}_{R}(M/\mathfrak{ a}M, \text{Hom}_{R}(X,Y))\cong \text{Hom}_{R}(M/\mathfrak{ a}M\otimes_{R}X, Y) 
$$
\cong \text{Hom}_{R}(M/\mathfrak{ a}^{2}M\otimes_{R}X, Y)\cong \text{Hom}_{R}(M/\mathfrak{ a}^{2}M, \text{Hom}_{R}(X,Y)).
$}
The conclusion $\text{Hom}_{R}(X,Y)\in\mathfrak{R}^{M}_{\mathfrak{ a}}$ follows by Proposition \ref{sec1:Prop- strongly I-coreduced}. For the second statement, the associative property of tensor products and Proposition \ref{sec1:Prop- strongly I-coreduced} are used, i.e., 
$M/\mathfrak{ a}M\otimes_{R}(X\otimes_{R}Y)\cong (M/\mathfrak{ a}M\otimes_{R} X)\otimes_{R} Y\cong (M/\mathfrak{ a}^{2}M\otimes_{R}X)\otimes_{R} Y\cong
M/\mathfrak{ a}^{2}M\otimes_{R}( X\otimes_{R}Y).$
So, $X\otimes_{R}Y\in \mathfrak{C}^{M}_{\mathfrak{ a}}$. 
\end{prf}

\section{Some closure properties}
 
\begin{prop} \label{Prop: Properties}
Let $\mathfrak{a}$ be an ideal of  $R$ and $M$ an $R$-module. Let $\{M_{j}\}_{j\in\Lambda}$ (resp. $\{N_{k}\}_{k\in\Lambda}$) be a family of $\mathfrak{ a}$-reduced modules with respect to $M$ (resp. $\mathfrak{ a}$-coreduced modules with respect to $M$) where $\Lambda$ is an index set. The following statements hold. 
\begin{itemize}
\item [$(1)$] {\normalfont $\prod_{j\in\Lambda}^{} M_{j} \in\mathfrak{R}^{M}_{\mathfrak{ a}}$} and  {\normalfont $\bigoplus_{k\in\Lambda} N_{k} \in\mathfrak{C}^{M}_{\mathfrak{ a}}$}. 
\item [$(2)$] If $\{M_{j}\}_{j\in\Lambda}$ is an inverse system of $R$-modules, then {\normalfont $\underset{j}\varprojlim M_{j}\in\mathfrak{R}^{M}_{\mathfrak{ a}}$}. 

\item[$(3)$] If $\{N_{k}\}_{k\in\Lambda}$ is a direct system of $R$-modules, then {\normalfont $\underset{k}\varinjlim N_{k}\in\mathfrak{C}^{M}_{\mathfrak{ a}}$}. 
\item [$(4)$] {\normalfont $\mathfrak{R}^{M}_{\mathfrak{ a}}$} (resp. {\normalfont $\mathfrak{C}^{M}_{\mathfrak{ a}}$}) is closed under submodules (resp. under quotients). 
  
\end{itemize}
\end{prop}
\begin{prf}
\begin{itemize}
\item [$(1)$] Since {\normalfont $M_{j}\in\mathfrak{R}^{M}_{\mathfrak{ a}}$} for each $j$, there is an isomorphism $\text{Hom}_{R}(M/\mathfrak{ a}M, M_{j})\cong \text{Hom}_{R}(M/\mathfrak{ a}^{2}M, M_{j})$ [Proposition \ref{Prop: reduced w.r.t M}]. This together with a property of the Hom functor that it commutes with arbitrary direct products in the second variable yields
$\text{Hom}_{R}(M/\mathfrak{ a}M, \prod_{j\in\Lambda} M_{j})\cong  \prod_{j\in\Lambda}\text{Hom}_{R}(M/\mathfrak{ a}M, M_{j})\cong \prod_{j\in\Lambda}\text{Hom}_{R}(M/\mathfrak{ a}^{2}M, M_{j})	
\cong \text{Hom}_{R}(M/\mathfrak{ a}^{2}M, \prod_{j\in\Lambda} M_{j}).$ This asserts the first statement. Suppose that $ N_{k}\in\mathfrak{C}^{M}_{\mathfrak{ a}}$ so that $M/\mathfrak{ a}M\otimes_{R} N_{k}\cong M/\mathfrak{ a}^{2}M\otimes_{R}N_{k}$ for each $k$ [Proposition \ref{sec1:Prop- strongly I-coreduced}]. Applying the tensor product property, that is, commutation with arbitrary direct sums:
$(M/\mathfrak{ a}M \otimes_{R}\bigoplus_{k\in\Lambda} N_{k})\cong \bigoplus_{k\in\Lambda} (M/\mathfrak{ a}M \otimes_{R} N_{k})\cong \bigoplus_{k } (M/\mathfrak{ a}^{2}M\otimes_{R} N_{k})	
\cong (M/\mathfrak{ a}^{2}M\otimes_{R} \bigoplus_{k\in\Lambda} N_{k}).$
\item [$(2)$] 
We use the fact that the functor $\text{Hom}_{R}(M/\mathfrak{ a}M,-)$ preserves inverse limits. So,
$\text{Hom}_{R}(M/\mathfrak{ a}M, \underset{j}\varprojlim M_{j})\cong \underset{j}\varprojlim  (\text{Hom}_{R}(M/\mathfrak{ a}M,M_{j}))\cong \underset{j}\varprojlim  (\text{Hom}_{R}(M/\mathfrak{ a}^{2}M,M_{j}))$
$\cong \text{Hom}_{R}(M/\mathfrak{ a}^{2}M,\underset{j}\varprojlim M_{j}).$
 
\item[$(3)$] Dually, since tensor products commute with direct limits:
$M/\mathfrak{ a}M\otimes_{R} \underset{k}\varinjlim  N_{k}\cong \underset{k}\varinjlim (M/\mathfrak{ a}M\otimes  N_{k})\cong \underset{k}\varinjlim (M/\mathfrak{ a}^{2}M\otimes  N_{k} ) 
\cong M/\mathfrak{ a}^{2}M\otimes_{R} \underset{k}\varinjlim  N_{k}.$
 
\item [$(4)$] By Definition \ref{Defn: reduced w.r.t M},  $\text{Hom}_{R}(M,N)$ is $\mathfrak{ a}$-reduced for each {\normalfont $N\in\mathfrak{R}^{M}_{\mathfrak{ a}}$}. Since for any submodule $X$ of $N$, the module $\text{Hom}_{R}(M,X)$ is a submodule of $\text{Hom}_{R}(M,N)$ and $\mathfrak{ a}$-reduced modules are closed under  submodules, the proof of the first statement follows. 
To prove the remaining statement, let $K$ be a submodule of $N$ and {\normalfont $N\in\mathfrak{C}^{M}_{\mathfrak{ a}}$}. Then $M\otimes_{R}N$ is $\mathfrak{ a}$-coreduced [Definition \ref{Defn: Cored w.r.t M}]. To complete the proof, it is enough to show that $M\otimes_{R}N/K$ is an $\mathfrak{ a}$-coreduced module. Since the functor $M\otimes_{R}-$ is right exact, applying it to a short exact sequence
$0\rightarrow K\rightarrow N\rightarrow N/K\rightarrow 0$, we get an exact sequence
$M\otimes_{R}K\rightarrow M\otimes_{R}N\rightarrow M\otimes_{R}N/K\rightarrow 0.$ It follows that $M\otimes_{R}N/K\cong (M\otimes_{R}N)/(M\otimes_{R}K)$,  which is a quotient of an $\mathfrak{ a}$-coreduced module $M\otimes_{R}N$. Since $\mathfrak{ a}$-coreduced modules are closed under quotients \cite[Proposition 5.1.5]{Annet's dissertation}, $M\otimes_{R}N/K$ is $\mathfrak{ a}$-coreduced and therefore $N/K$ is $\mathfrak{ a}$-coreduced with respect to $M$.  
\end{itemize}	
\end{prf}

\begin{paragraph}\noi Recall that a category $\mathcal{C}$ is \textit{complete} if $\underset{}\varprojlim M_{j}$ exists in $\mathcal{C}$ for every inverse system $\{M_{j}\}_{j\in\Lambda}$ in $\mathcal{C}$. Dually,  $\mathcal{C}$ is \textit{cocomplete} if $\underset{}\varinjlim N_{j}$ exists in $\mathcal{C}$ for every direct system $\{N_{j}\}_{j\in\Lambda}$ in $\mathcal{C}$; see \cite[Definition 5.27]{J.Rotman- Introduction hom.algebra}. Moreover, a class $\mathscr{C}$ of objects is called a \textit{pretorsion class} if it is closed under quotient objects and coproducts, and is a \textit{pretorsion-free class} if it is closed under subobjects and products \cite[page 137]{Stenstrom--rings of quotients}. By Proposition \ref{Prop: Properties} it follows that
\end{paragraph}
\begin{prop}\label{prop: complete vs cocomplete}
The full subcategory  {\normalfont $\mathfrak{R}^{M}_{\mathfrak{ a}}$} (resp. {\normalfont $\mathfrak{C}^{M}_{\mathfrak{ a}}$}) of {\normalfont $R\text{-Mod}$} is complete and forms a pretorsion-free class (resp. cocomplete and forms a pretorsion class). 
\end{prop}
\begin{lem}{\normalfont\cite[Proposition 2.7]{Annet-David : generalised reduced}} \label{Lem: localizedreduced modules}
Let $\mathfrak{a}$ be an ideal of  $R$, and let $S$ be a multiplicative subset of $R$. An $R$-module $N$ is $\mathfrak{a}$-reduced (resp. reduced) if and only if an  $S^{-1}R$-module $S^{-1}N$ is  $S^{-1}\mathfrak{a}$-reduced (resp. reduced).
\end{lem}
\begin{paragraph}\noi Statements in Lemma \ref{lem: Properties of Localisation} are well-known.
\end{paragraph}
\begin{lem}\label{lem: Properties of Localisation}
Let $S$ be a multiplicative subset of $R$, and let $M$ and $N$ be $R$-modules.
\begin{itemize} 
\item[$(1)$] If $M$ is finitely presented, then 
{\normalfont $S^{-1}(\text{Hom}_{R}(M,N))\cong \text{Hom}_{S^{-1}R}(S^{-1}M, S^{-1}N).$
\item[$(2)$] $ S^{-1}(M\otimes_{R}N)\cong S^{-1}M \otimes_{S^{-1}R} S^{-1}N.$} 
\item[$(3)$]  $S^{-1}(R/\mathfrak{a})\cong S^{-1}R/S^{-1}\mathfrak{a}.$ 
\item[$(4)$] $S^{-1}\Lambda_{\mathfrak{a}}(N)\cong \Lambda_{S^{-1}\mathfrak{a}}(S^{-1}N).$ 
\end{itemize} 
\end{lem}
\begin{lem}\label{Lem: Localization of coreduced modules}
Let $\mathfrak{a}$ be an ideal of  $R$. An $R$-module $N$ is $\mathfrak{a}$-coreduced if and only if $S^{-1}N$ over $S^{-1}R$ is $S^{-1}\mathfrak{a}$-coreduced. 
\end{lem}
\begin{prf}
Recall that $N\in\mathfrak{ C}_{\mathfrak{ a}}$ if and only if $\Lambda_{\mathfrak{a}}(N)\cong R/\mathfrak{a}\otimes_{R}N$ \cite[Proposition 2.3]{David-Application-I}. This together with Lemma \ref{lem: Properties of Localisation}, $\Lambda_{S^{-1}\mathfrak{a}}(S^{-1}(N))\cong S^{-1}(\Lambda_{\mathfrak{a}}(N))\cong S^{-1}(R/\mathfrak{a}\otimes_{R}N)\cong S^{-1}(R/\mathfrak{a})\otimes_{S^{-1}R} S^{-1}N \cong S^{-1}R/S^{-1}\mathfrak{a}\otimes_{S^{-1}R} S^{-1}N$ if and only if $S^{-1}N\in\mathfrak{C}_{S^{-1}\mathfrak{ a}}$.
\end{prf}
\begin{paragraph}\noi
It is worth mentioning that not all properties of $R$-modules localize. For instance, for an arbitrary commutative ring  $R$ the property of being an injective $R$-module does not localize \cite[Section A.5]{Peter-Schenzez}.
\end{paragraph}
\begin{prop}\label{Prop: Localised modules}
Let $S$ be a multiplicatively closed subset of $R$ and {\normalfont $M,N\in R\text{-Mod}$}.
\begin{itemize} 
\item[$(1)$] If $M$ is finitely presented, then $N\in\mathfrak{R}^{M}_{\mathfrak{ a}}$ if and only if $S^{-1}N\in\mathfrak{R}^{S^{-1}M}_{S^{-1}\mathfrak{ a}}$.
\item[$(2)$] $N\in\mathfrak{C}^{M}_{\mathfrak{ a}}$ if and only if $S^{-1}N\in\mathfrak{C}^{S^{-1}M}_{S^{-1}\mathfrak{ a}}$.
\end{itemize}
\end{prop}
\begin{prf}
\begin{itemize}
\item [$(1)$] By definition, $N\in\mathfrak{R}^{M}_{\mathfrak{ a}}$  if and only if $\text{Hom}_{R}(M,N)\in\mathfrak{R}^{}_{\mathfrak{ a}}$ if and only if $S^{-1}\text{Hom}_{R}(M,N)\in\mathfrak{R}_{S^{-1}\mathfrak{ a}}$ [Lemma \ref{Lem: localizedreduced modules}]. The conclusion is by Lemma \ref{lem: Properties of Localisation} $(1)$.
\item[$(2)$]
 
$N\in\mathfrak{C}^{M}_{\mathfrak{ a}}$  if and only if $M\otimes_{R}N\in\mathfrak{C}_{\mathfrak{ a}}$  if and only if $S^{-1}(M\otimes_{R}N)\cong S^{-1}M\otimes_{S^{-1}R}S^{-1}N\in\mathfrak{C}_{S^{-1}\mathfrak{ a}}$ by Lemma \ref {lem: Properties of Localisation} $(2)$ and Lemma \ref{Lem: Localization of coreduced modules}. Thus, $S^{-1}N\in \mathfrak{C}^{S^{-1}M}_{S^{-1}\mathfrak{ a}}$. 
\end{itemize}
\end{prf}

\begin{paragraph}\noi
Neither $\mathfrak{R}^{M}_{\mathfrak{ a}}$ nor $\mathfrak{C}^{M}_{\mathfrak{ a}}$ is closed under extension, in general. 
	First recall that $\mathfrak{R}_{\mathfrak{ a}}\subseteq\mathfrak{R}^{M}_{\mathfrak{ a}}$, and $\mathfrak{C}_{\mathfrak{ a}}\subseteq\mathfrak{C}^{M}_{\mathfrak{ a}}$. Given the sequence of $\Z$-modules $0\rightarrow (2)\Z_{4}\rightarrow \Z_{4}\rightarrow \Z_{4}/(2)\Z_{4}\rightarrow 0$, we have $(2)\Z_{4}, \Z_{4}/(2)\Z_{4}\in\mathfrak{R}^{\Z}_{(2)}$, but $\Z_{4}\notin\mathfrak{R}^{\Z}_{(2)}$. Thus, $\mathfrak{R}^{M}_{\mathfrak{ a}}$ is not closed under extension. For the remaining case, in the exact sequence $0\rightarrow 2\Z/4\Z\rightarrow \Z/4\Z\rightarrow \Z/2\Z\rightarrow 0$ of $\Z$-modules, $2\Z/4\Z,\Z/2\Z\in\mathfrak{C}^{\Z}_{(2)}$. However, $\Z/4\Z\notin\mathfrak{ C}^{\Z}_{(2)}$. Thus, $\mathfrak{C}^{M}_{\mathfrak{ a}}$ is not closed under extension.
	\end{paragraph}

\section{Some Duality Theorems}
\begin{prop}\label{Prop: leading to duality}
	Let $\mathfrak{ a}$ be an ideal of  $R$, and let $X,Y\in{\normalfont R\text{-Mod}}$ where $Y$ is an injective cogenerator of {\normalfont $R\text{-Mod}$}.
	\begin{itemize}
		\item [$(1)$] The $R$-module {\normalfont $X\in  \mathfrak{C}^{M}_{\mathfrak{ a}}$} if and only if {\normalfont $\text{Hom}_{R}(X,Y)\in\mathfrak{R}^{M}_{\mathfrak{ a}}$}. 
		\item[$(2)$] Let $\mathfrak{ a}$ be finitely generated  and $M$ a finitely generated $R$-module. The $R$-module  {\normalfont $\text{Hom}_{R}(X,Y)\in\mathfrak{C}^{M}_{\mathfrak{ a}}$} for all {\normalfont $X\in\mathfrak{R}^{M}_{\mathfrak{ a}}$}.
	\end{itemize} 
\end{prop}
\begin{prf}
	\begin{itemize}
		\item[$(1)$] The necessity is proved in Proposition \ref{Hom and tensor is str. (co)reduced}. We prove only  the sufficiency. Suppose that {\normalfont $\text{Hom}_{R}(X,Y)\in\mathfrak{R}^{M}_{\mathfrak{ a}}$}. This together with the Hom-Tensor adjunction gives:
		\begin{equation*}
			{\normalfont \text{Hom}_{R}(M/\mathfrak{ a}M\otimes_{R}X,Y)\cong  \text{Hom}_{R}(M/\mathfrak{ a}M,\text{Hom}_{R}(X,Y))\cong }	
		\end{equation*}
		\begin{equation*}
			\text{Hom}_{R}(M/\mathfrak{ a}^{2}M,\text{Hom}_{R}(X,Y))\cong \text{Hom}_{R}(M/\mathfrak{ a}^{2}M\otimes_{R}X,Y). 
		\end{equation*} 
		But since $Y$ is an injective cogenerator the functor $\text{Hom}_{R}(-,Y)$ reflects isomorphisms. So, $ M/\mathfrak{ a}M\otimes_{R}X \cong M/\mathfrak{ a}^{2}M\otimes_{R}X$. It follows by Proposition \ref{sec1:Prop- strongly I-coreduced} that $X\in\mathfrak{C}^{M}_{\mathfrak{ a}}$.
		\item[$(2)$]
		Since $M$ is finitely generated so is $M/\mathfrak{ a}^{k}M$ for every $k\in\Z^{+}$.	By virtue of \cite[11.1.3]{Peter-Schenzez} and by Proposition \ref{Prop: reduced w.r.t M}, we get $M/\mathfrak{ a}M\otimes_{R}\text{Hom}_{R}(X,Y)\cong \text{Hom}_{R}(\text{Hom}_{R}(M/\mathfrak{ a}M,X), Y)\cong$
		$ \text{Hom}_{R}(\text{Hom}_{R}(M/\mathfrak{ a}^{2}M,X), Y)\cong M/\mathfrak{ a}^{2}M \otimes_{R} \text{Hom}_{R}(X,Y).$ The conclusion $\text{Hom}_{R}(X,Y)\in\mathfrak{C}^{M}_{\mathfrak{ a}}$ holds by Proposition \ref{sec1:Prop- strongly I-coreduced}.
	\end{itemize}
\end{prf}
\begin{paragraph}\noi
	Denote by $N^{\vee}:= \text{Hom}_{R}(N,E)$ the general Matlis dual of an $R$-module $N$ where $E$ is an injective cogenerator of $R\text{-Mod}$. 
\end{paragraph}
\begin{cor}\label{Cor: X is cor iff its dual is red}
	Let $\mathfrak{ a}$ be an ideal of  $R$ and $X$ an $R$-module. Then {\normalfont $X\in\mathfrak{C}^{M}_{\mathfrak{ a}}$}  if and only if {\normalfont $X^{\vee}\in\mathfrak{R}^{M}_{\mathfrak{ a}}$}. If {\normalfont $X\in\mathfrak{R}^{M}_{\mathfrak{ a}}$}, then {\normalfont $X^{\vee}\in\mathfrak{C}^{M}_{\mathfrak{ a}}$} provided $\mathfrak{ a}$ is finitely generated and $M$ is a finitely generated $R$-module. 
	 
\end{cor}
\begin{prf}	By taking $Y:= E$, Proposition \ref{Prop: leading to duality} $(1)$ (resp. $(2)$) provides immediate proof of the first (resp. the second) statement.	
\end{prf}

\begin{thm}\label{Prop: Duality_ Gen Gamma vs Gen. Lambda }
	Suppose that $\mathfrak{ a}$ and  $M$  are finitely generated. If {\normalfont $N\in\mathfrak{R}^{M}_{\mathfrak{a}}$}, then $$\Gamma_{ \mathfrak{a}}(M,N)^{\vee}\cong \Lambda_{ \mathfrak{a}}(M,N^{\vee}).$$
\end{thm}
\begin{prf}
	By the isomorphism evaluation,  $\text{Hom}_{R}(\text{Hom}_{R}(M/\mathfrak{a}M,N),E)\cong M/\mathfrak{ a}M\otimes_{R}\text{Hom}_{R}(N,E)$ \cite[11.1.3]{Peter-Schenzez}. Since $N\in\mathfrak{R}^{M}_{\mathfrak{ a}}$, $\Gamma_{ \mathfrak{a}}(M,N)\cong \text{Hom}_{R}(M/\mathfrak{ a}M,N)$ [Proposition \ref{Prop: reduced w.r.t M}], and by Corollary \ref{Cor: X is cor iff its dual is red}, $\Lambda_{ \mathfrak{a}} (M, N^{\vee})\cong M/\mathfrak{ a}M\otimes_{R}N^{\vee}$.  So,
	$\Gamma_{ \mathfrak{a}}(M,N)^{\vee}\cong (\text{Hom}_{R}(M/\mathfrak{ a}M,N))^{\vee}\cong \text{Hom}_{R}(\text{Hom}_{R}(M/\mathfrak{ a}M, N),E)$
	$\cong M/\mathfrak{ a}M\otimes_{R}\text{Hom}_{R}(N,E)\cong M/\mathfrak{ a}M\otimes_{R}N^{\vee}\cong \Lambda_{ \mathfrak{a}}(M,N^{\vee}).$  
\end{prf}

\begin{paragraph}\noi Let $i\in\Z$. We recall the Ext-Tor Duality \cite[1.4.1]{Peter-Schenzez}. There is an isomorphism  {\normalfont $\text{Hom}_{R}(\text{Tor}^{R}_{i}(M,N),J)\cong \text{Ext}^{i}_{R}(M, \text{Hom}_{R}(N,J))$} for all $i$ and all $R$-modules $M,N,J$ with $J$ injective. For a finitely generated ideal $\mathfrak{ a}$, $\Gamma_{ \mathfrak{a}}(N)^{\vee}\cong \Lambda_{ \mathfrak{a}}(N^{\vee})$ \cite[Observations 2.2.13]{Peter-Schenzez}. 
\end{paragraph}
\begin{thm} \label{Prop: Gen. Matlis Duality of generalised functors}
	Let $\mathfrak{ a}$ be an ideal of a ring $R$, and let $M$ and $N$ be $R$-modules. If {\normalfont $N\in\mathfrak{C}^{M}_{\mathfrak{a}}$}, then $$\Lambda_{ \mathfrak{a}}(M,N)^{\vee}\cong \Gamma_{ \mathfrak{a}}(M,N^{\vee}).$$
\end{thm}
\begin{prf}
	$(M\otimes_{R}N)^{\vee}\cong \text{Hom}_{R}(M,N^{\vee})$ is evident from the Ext-Tor Duality with $i=0$.
	Since $N\in\mathfrak{C}^{M}_{\mathfrak{a}}$, and thus {\normalfont $N^{\vee}\in\mathfrak{R}^{M}_{\mathfrak{ a}}$} by [Corollary \ref{Cor: X is cor iff its dual is red}], we get
	$\Lambda_{ \mathfrak{a}}(M,N)^{\vee}\cong (M/\mathfrak{ a}M\otimes_{R}N)^{\vee}\cong \text{Hom}_{R}(M/\mathfrak{ a}M,N^{\vee})\cong \Gamma_{ \mathfrak{a}}(M,N^{\vee})$ [Propositions \ref{Prop: reduced w.r.t M} and \ref{sec1:Prop- strongly I-coreduced}].
\end{prf}
\begin{paragraph}\noi Recall that an $R$-module $N$ is \textit{reflexive} if $N\cong N^{\vee\vee}$.
	\end{paragraph}
\begin{cor}
	Let $\mathfrak{ a}$ and $M$ be finitely generated. If $N\in\mathfrak{R}^{M}_{\mathfrak{ a}}\cap \mathfrak{ C}^{M}_{\mathfrak{ a}}$ is reflexive, then $\Gamma_{ \mathfrak{a}}(M,N)$ and $\Lambda_{ \mathfrak{a}}(M,N)$ are also reflexive $R$-modules.
\end{cor}
\begin{prf}
	The proof follows by Theorem \ref{Prop: Duality_ Gen Gamma vs Gen. Lambda } and Theorem \ref{Prop: Gen. Matlis Duality of generalised functors}.
\end{prf}

\begin{paragraph}\noi
	Let $M\in R\text{-Mod}$. Recall that {\normalfont $\mathfrak{R}_{\mathfrak{ a}}^{M}$} (resp. {\normalfont $\mathfrak{C}_{\mathfrak{ a}}^{M}$}) denotes all $\mathfrak{ a}$-reduced (resp. all $\mathfrak{ a}$-coreduced) $R$-modules with respect to $M$. 
	Denote $ \mathfrak{R}^{M}_{\mathfrak{ a}}\cap \mathfrak{C}^{M}_{\mathfrak{ a}}$ by $\mathfrak{B}^{M}_{\mathfrak{ a}}$.
\end{paragraph}
\begin{lem}\label{Lem: Modules in the intersection}
	Let $\mathfrak{ a}$ be an ideal of  $R$, and let $M,N$ be $R$-modules. Then the $R$-modules {\normalfont $M/\mathfrak{ a}M\otimes_{R}N, \text{Hom}_{R}(M/\mathfrak{ a}M,N)\in \mathfrak{B}_{\mathfrak{ a}}^{M}$}.
\end{lem}
\begin{prf}
	For any $M,N\in R\text{-Mod}$ and any ideal $\mathfrak{ a}$ of  $R$, the modules $M/\mathfrak{ a}M\otimes_{R}N\cong R/\mathfrak{ a}\otimes_{R}(M\otimes_{R}N)$ and $\text{Hom}_{R}(M/\mathfrak{ a}M,N)\cong \text{Hom}_{R}(R/\mathfrak{ a},\text{Hom}_{R}(M,N))$  are both $\mathfrak{ a}$-reduced and $\mathfrak{ a}$-coreduced \cite[Proposition 2.9]{David-Application-I}. It follows, by both Lemma \ref{Lem: further example-red.w.r.t M} and Proposition \ref{Lem: cor tensor with cor is cor}, that $M/\mathfrak{ a}M\otimes_{R}N$ and $\text{Hom}_{R}(M/\mathfrak{ a}M,N)$ are $\mathfrak{ a}$-reduced and $\mathfrak{ a}$-coreduced with respect to any module. 
\end{prf}
\begin{cor}\label{Cor: Hom and tensor are both red. and cored} Let $\mathfrak{ a}$ be an ideal of $R$. For all $N\in\mathfrak{R}^{M}_{\mathfrak{ a}}$ and all $P\in\mathfrak{C}^{M}_{\mathfrak{ a}}$, $\Gamma_{ \mathfrak{a}}(M,N)\in\mathfrak{C}^{M}_{\mathfrak{ a}}$, and $\Lambda_{ \mathfrak{a}}(M,P)\in\mathfrak{R}^{M}_{\mathfrak{ a}}$.	
\end{cor}
\begin{prf}
	This holds by Lemma \ref{Lem: Modules in the intersection}, since for any $N\in\mathfrak{R}^{M}_{\mathfrak{ a}}$ and any $P\in\mathfrak{C}^{M}_{\mathfrak{ a}}$, $\Gamma_{ \mathfrak{a}}(M,N)\cong\text{Hom}_{R}(M/\mathfrak{a}M,N)$, and $\Lambda_{ \mathfrak{a}}(M,P)\cong M/\mathfrak{ a}M\otimes_{R}P$ [Propositions \ref{Prop: reduced w.r.t M} and \ref{sec1:Prop- strongly I-coreduced}].
\end{prf}
\begin{paragraph}\noi Let $M,N\in R\text{-Mod}$. The functors $\text{Hom}_{R}(M,-)$ and $M\otimes-$ are adjoint. The naive expectation is that their generalisations $\Gamma_{ \mathfrak{a}}$ (resp.  $\Lambda_{ \mathfrak{a}}$) are also adjoint. Unfortunately, this is not the case in general. The conditions under which they become adjoint has been a subject of interest in different contexts; see for instance \cite{Alonso-Lipman} in the context of sheaves, \cite{porta2014homology, Peter-Schenzez} in the context of derived categories over modules, and \cite{David-Application-I} in the context of modules. In each of these cases, the adjunction is called the \textit{Greenlees-May Duality} (\textit{GM-Duality} for short). In Theorem \ref{GGM-Duality} below, we demonstrate that this Duality also holds in the setting of $\mathfrak{R}^{M}_{\mathfrak{ a}}$ and $\mathfrak{C}^{M}_{\mathfrak{ a}}$ which is an extension of \cite[Theorem 3.4]{David-Application-I} proved in the setting of $\mathfrak{R}_{\mathfrak{ a}}$ and $\mathfrak{ C}_{\mathfrak{ a}}$.
	\end{paragraph}
\begin{thm}{\normalfont\textbf{ [The Generalised Greenlees-May Duality]}}\label{GGM-Duality}
	Let $M$ be an $R$-module. Define {\normalfont $\Gamma_{ \mathfrak{a}}(M,-): \mathfrak{R}_{\mathfrak{ a}}^{M}\rightarrow \mathfrak{C}_{\mathfrak{ a}}^{M}$} by $N\mapsto \Gamma_{ \mathfrak{a}}(M,N)$ and {\normalfont $\Lambda_{ \mathfrak{a}}(M,-): \mathfrak{C}_{\mathfrak{ a}}^{M}\rightarrow \mathfrak{R}_{\mathfrak{ a}}^{M}$} by $P\mapsto \Lambda_{ \mathfrak{a}}(M,P)$. Then for all {\normalfont $N\in \mathfrak{R}_{\mathfrak{ a}}^{M}$} and all {\normalfont $P\in\mathfrak{C}_{\mathfrak{ a}}^{M}$}, {\normalfont
		$$\text{Hom}_{R}(\Lambda_{ \mathfrak{a}}(M, P),N)\cong \text{Hom}_{R}(P, \Gamma_{ \mathfrak{a}}(M,N)).$$}
\end{thm}
\begin{prf} By Corollary \ref{Cor: Hom and tensor are both red. and cored},
	$\Gamma_{ \mathfrak{a}}(M,N)\cong \text{Hom}_{R}(M/\mathfrak{ a}M,N)\in \mathfrak{C}_{\mathfrak{ a}}^{M}$ and $\Lambda_{ \mathfrak{a}}(M,P)\cong M/\mathfrak{ a}M\otimes_{R}P\in \mathfrak{R}_{\mathfrak{ a}}^{M}$ for all $N\in \mathfrak{R}_{\mathfrak{ a}}^{M}$ and all $P\in \mathfrak{C}_{\mathfrak{ a}}^{M}$. Since $M/\mathfrak{ a}M\otimes_{R}P\cong P\otimes_{R}M/\mathfrak{ a}M$, there are isomorphisms
	$
	\text{Hom}_{R}(\Lambda_{ \mathfrak{a}}(M,P),N)\cong \text{Hom}_{R}(M/\mathfrak{ a}M\otimes_{R}P,N)\cong \text{Hom}_{R}(P\otimes_{R}M/\mathfrak{ a}M,N)
	\cong \text{Hom}_{R}(P,\text{Hom}_{R}(M/\mathfrak{ a}M,N))\cong \text{Hom}_{R}(P, \Gamma_{ \mathfrak{a}}(M,N)).$
	\end{prf}
	\begin{paragraph}\noi 
		The functor $\Lambda_{ \mathfrak{a}}(M,-):R \text{-Mod}\rightarrow R \text{-Mod}$ which is given by $N\mapsto \Lambda_{ \mathfrak{a}}(M,N):=\underset{k} \varprojlim (M/\mathfrak{ a}^{k}M \otimes_{R}N)$ is not exact; it is neither left exact nor right exact, since it is the composition of the left exact functor $\varprojlim$, and the right exact functor $M\otimes_{R}-$, see \cite{T.T. Nam-generalised I-adic completion}.
	\end{paragraph}
	\begin{cor}\label{Lem: Exactness of the functors}
		{\normalfont Let $M$ be an $R$-module.
			\begin{itemize}
				\item [$(1)$] The functor $\Gamma_{ \mathfrak{a}}(M,-): \mathfrak{R}_{\mathfrak{a}}^{M}\rightarrow\mathfrak{C}_{\mathfrak{ a}}^{M}$ is left exact. 
				\item [$(2)$] The functor $\Lambda_{ \mathfrak{a}}(M,-): \mathfrak{C}_{\mathfrak{ a}}^{M}\rightarrow \mathfrak{R}_{\mathfrak{ a}}^{M}$ is right exact. 
		\end{itemize}}
	\end{cor}
	\begin{prf}
		This immediately follows from Theorem \ref{GGM-Duality} and \cite[Theorem 2.6.1]{weibel1995introduction}.
	\end{prf}
	\begin{paragraph}\noi
	Denote by $\textbf{Set}$ the category of all sets. Recall that a functor $F: \mathscr{C}\rightarrow \textbf{Set}$ is \textit{representable} if $F$ is naturally isomorphic to $\text{Hom}_{R}(K,-)$ for some $K\in \mathscr{C}$.
	\end{paragraph}
	\begin{prop}\label{prop: representablity of gen. Gamma functor}
		Let $\mathfrak{ a}$ be an ideal of  $R$, and let $M$ be an $R$-module. The functor {\normalfont $\Gamma_{ \mathfrak{a}}(M,-):\mathfrak{R}_{\mathfrak{ a}}^{M}\rightarrow\mathfrak{C}_{\mathfrak{ a}}^{M}$} is representable, and therefore preserves all limits. In particular, {\normalfont $\Gamma_{ \mathfrak{a}}(M,-)$} preserves inverse limits, coproducts, pullbacks, equalizers, and  terminal objects. 
	\end{prop}
	\begin{prf}
		Since, by Lemma \ref{Lem: further example-red.w.r.t M},  $M/\mathfrak{a}M\in \mathfrak{R}_{\mathfrak{ a}}^{M}$ and  $\Gamma_{ \mathfrak{a}}(M,N)\cong \text{Hom}_{R}(M/\mathfrak{ a}M,N)$ [Proposition \ref{Prop: reduced w.r.t M}], $\Gamma_{ \mathfrak{a}}(M,-)$ is representable. Moreover, since $\Gamma_{ \mathfrak{a}}(M,-)$ is right adjoint to $\Lambda_{ \mathfrak{a}}(M,-)$ by Theorem \ref{GGM-Duality}, it preserves limits; see \cite[Theorem 2.6.10]{weibel1995introduction}.
	\end{prf}
	
	\begin{paragraph}\noi
		We denote by $\text{Nat}(-,-)$ the set of all natural transformations between any two functors.
	\end{paragraph}	
	\begin{prop}\label{Prop: 4 natural transformation}
		Let $\mathfrak{ a}$ be an ideal of  $R$, and let $M$ be an $R$-module. For every functor {\normalfont $\Gamma_{ \mathfrak{a}}(M,-):\mathfrak{R}_{\mathfrak{ a}}^{M}\rightarrow \mathfrak{C}_{\mathfrak{ a}}^{M}$}, there is an isomorphism {\normalfont $\text{Nat}(\Gamma_{ \mathfrak{a}}(M,-), \Gamma_{ \mathfrak{a}}(M,-))\cong \text{End}_{R}(M/\mathfrak{ a}M).$}
	\end{prop}
	\begin{prf}
		Consider the functor $\text{Hom}_{R}(M/\mathfrak{ a}M,-): \mathfrak{R}_{\mathfrak{a}}^{M}\rightarrow \textbf{Set}$.  Since $\Gamma_{ \mathfrak{a}}(M,-)$ is representable, by Yoneda's Lemma there is a one-to-one correspondence between 
		$\text{Nat}(\text{Hom}_{R}(M/\mathfrak{ a}M,-),\Gamma_{ \mathfrak{a}}(M,-) )$ and $\Gamma_{ \mathfrak{a}}(M,M/\mathfrak{ a}M)$. This together with the isomorphism $\Gamma_{ \mathfrak{a}}(M,-)\cong \text{Hom}_{R}(M/\mathfrak{ a}M,-)$ yield
		$\text{Nat}(\Gamma_{ \mathfrak{a}}(M,-), \Gamma_{ \mathfrak{a}}(M,-))\cong \Gamma_{ \mathfrak{a}}(M,M/\mathfrak{ a}M)\cong \text{Hom}_{R}(M/\mathfrak{ a}M,M/\mathfrak{ a}M)$$=\text{End}_{R}(M/\mathfrak{ a}M).$
	\end{prf}
	\begin{prop}
		The functor $\Lambda_{ \mathfrak{a}}(M,-): \mathfrak{ C}^{M}_{\mathfrak{ a}}\rightarrow \mathfrak{R}^{M}_{\mathfrak{ a}}$ preserves  all colimits (in particular, coproducts, direct limits, cokernels, pushouts, initial objects and coequilizers). 
	\end{prop}
	\begin{prf}
		By 	Corollary \ref{Lem: Exactness of the functors}, $\Lambda_{ \mathfrak{a}}(M,-)$ is left adjoint to $\Gamma_{ \mathfrak{a}}(M,-)$. So, by \cite[Theorem 2.6.10]{weibel1995introduction}, it preserves colimits.
	\end{prf}

\section{Effect on generalised local (co)homology}

\begin{paragraph}\noi
 The purpose of this section is to give some applications of (co)reduced modules with respect to a given module $M$ on generalised local cohomology (resp. generalised local homology) modules.  In this section, all rings $R$ are Noetherian.
\end{paragraph}

\begin{paragraph}\noi
Recall that $\Gamma_{ \mathfrak{a}}(M,-)$ is a left exact endo-functor on $R\text{-Mod}$ [Lemma \ref{Isomorphisms: The Gammas and the Lamdas}]; see also \cite{T.T. Nam-generalised I-adic completion}. The right derived functor of $\Gamma_{ \mathfrak{a}}(M,-)$ is denoted by $\text{H}^{i}_{\mathfrak{ a}}(M,-)$ and called the \textit{generalised local cohomology functor with respect to $\mathfrak{ a}$}. The \textit{ $i^{\text{th}}$-generalised local cohomology} module $\text{H}^{i}_{\mathfrak{a}}(M,N)$ of $R$-modules $M,N$ with respect to $\mathfrak{a}$ which was introduced by Herzog \cite{Herzog-generalised cohomology} is defined by $\text{H}^{i}_{\mathfrak{a}}(M,N):= \underset{k}\varinjlim~ \text{Ext}^{i}_{R}(M/\mathfrak{a}^{k}M, N).$ 
Dually, the {\it $i^{\text{th}}$-generalised local homology} module of $M,N$  with respect to $\mathfrak{ a}$ is given by 
{\normalfont $\text{H}^{\mathfrak{ a}}_{i}(M,N):=\underset{k}\varprojlim~ \text{Tor}^{R}_{i}(M/\mathfrak{a}^{k}M, N)$} \cite[page 446]{T.T. Nam-generalised I-adic completion}.
Note that if $M= R$ then the generalised local cohomology coincides with the Grothendieck local cohomology which is given by   $ \text{H}^{i}_{\mathfrak{ a}}(N):= \underset{k}\varinjlim\text{ Ext}^{i}_{R}(R/\mathfrak{ a}^{k},N)$ for any $R$-module $N$. Similarly, when $M=R$ the generalised local homology coincides with the local homology defined by $H_{i}^{\mathfrak{ a}}(N)\cong \underset{k}\varprojlim\text{ Tor}^{R}_{i}(R/\mathfrak{ a}^{k},N)$ for all $R$-modules $N$. For the basics of local (co)homology modules we refer to \cite{brodmann2013local, Hartshorne-Local cohomology, Peter-Schenzez} among others.  
\end{paragraph}
 
\begin{paragraph}\noi
If $M$ is a finitely generated $R$-module and $N$ is an $\mathfrak{ a}$-torsion $R$-module, then $\text{H}^{i}_{\mathfrak{ a}}(M,N)\cong \text{Ext}^{i}_{R}(M,N)$ for all $i\ge 0$ \cite[Lemma 1.1]{Yassemi S- Associated primes}. For $R$-modules $\mathfrak{ a}$-reduced with respect to some module $M$, we have Proposition \ref{Prop: generalised local cohomology}.
\end{paragraph}
\begin{prop}\label{Prop: generalised local cohomology}
Let $\mathfrak{a}$ be an ideal of  $R$, and let $M$ be an $R$-module. If $N\in\mathfrak{R}^{M}_{\mathfrak{ a}}$, then  
{\normalfont $\text{H}^{i}_{\mathfrak{a}}(M,N)\cong \text{Ext}^{i}_{R}(M/\mathfrak{a}M, N)$ for all $i\ge0$}.
\end{prop}
\begin{prf}
Since $N\in\mathfrak{R}_{\mathfrak{ a}}^{M}$, there is an isomorphism $\Gamma_{\mathfrak{ a}}(M,N)\cong \text{Hom}_{R}(M/\mathfrak{a}M, N)$ [Proposition \ref{Prop: reduced w.r.t M}]. It follows that the right derived functor of $\Gamma_{\mathfrak{ a}}(M,N)$ which is $\text{H}^{i}(M,N)$ will be given by $\text{H}^{i}_{\mathfrak{ a}}(M,N)
\cong \text{Ext}^{i}_{R}(M/\mathfrak{a}M, N)$ for all $i\ge 0$. 
\end{prf}
\begin{cor}\label{cor: to generalised local cohomology-vanishing of generalised local cohomology}
Let $M$ and $N$ be $R$-modules.	
If $M/\mathfrak{a}M$ is a projective $R$-module, then {\normalfont $\text{H}^{i}_{\mathfrak{a}}(M,N)=0$ for all $N\in\mathfrak{R}^{M}_{\mathfrak{ a}}$ and all $i\ge 1$.}
\end{cor} 
\begin{prf}
If the $R$-module $ M/\mathfrak{ a}M$ is projective, then  by the general theory of the Ext-functor, $\text{Ext}^{i}_{R}(M/\mathfrak{ a}M, N)=0$ for all $i\ge 1$. Since $\text{H}^{i}_{\mathfrak{ a}}(M,N)\cong \text{Ext}^{i}_{R}(M/\mathfrak{ a}M, N)$ for all $i\ge0$ by Proposition \ref{Prop: generalised local cohomology}, the result follows.
\end{prf}  
\begin{paragraph}\noi For a finitely generated $R$-module $M$ and an $\mathfrak{ a}$-adically complete artinian module $N$ over a local ring $(R,\mathfrak{ m})$, there is an isomorphism $\text{H}^{\mathfrak{a}}_{i}(M,N)\cong \text{Tor}_{i}^{R}(M,N)$ \cite[Lemma 2.7]{T T Nam-GLH for artinian modules}. For modules in $\mathfrak{ C}^{M}_{\mathfrak{ a}}$ we have Proposition \ref{Prop: generalised local homology}.
\end{paragraph}
\begin{prop}\label{Prop: generalised local homology}
Let $\mathfrak{a}$ be an ideal of  $R$, $M$ any $R$-module, and let $N\in\mathfrak{ C}^{M}_{\mathfrak{ a}}$. Then
{\normalfont $\text{H}^{\mathfrak{ a}}_{i}(M,N)\cong \text{Tor}^{R}_{i}(M/\mathfrak{a}M, N)$ for all $i\ge0$}.
\end{prop}
\begin{prf}
By Proposition \ref{sec1:Prop- strongly I-coreduced}, $\Lambda_{\mathfrak{ a}}(M,N)\cong M/\mathfrak{ a}M\otimes_{R}N.$ It follows that 
$\text{H}^{\mathfrak{ a}}_{i}(M,N)
\cong \text{Tor}^{R}_{i}(M/\mathfrak{a}M, N) $ for all $i$.
\end{prf}

\begin{cor}\label{Cor: Vanishing of gen.local homology over coreduced}
Let $\mathfrak{ a}$ be an ideal of a ring $R$, and {\normalfont $M,N\in R\text{-Mod}$}. If $M/\mathfrak{a}M$ is flat and $N\in\mathfrak{ C}^{M}_{\mathfrak{ a}}$, or $N$ is flat and $N\in\mathfrak{ C}^{M}_{\mathfrak{ a}}$,  then {\normalfont $\text{H}^{\mathfrak{ a}}_{i}(M,N)=0$ for all $i>0$.}
 
\end{cor}
\begin{prf} This is immediate by Proposition \ref{Prop: generalised local homology}, since $\text{H}^{\mathfrak{ a}}_{i}(M,N)\cong \text{Tor}^{R}_{i}(M/\mathfrak{ a}M,N)$ and $\text{Tor}^{R}_{i}(M,N)=0$ for all $i>0$ if either $M$ or $N$ is flat \cite[Lemma 3.2.8]{weibel1995introduction}.
\end{prf}

\begin{paragraph}\noi Propositions \ref{Prop: generalised local cohomology} and \ref{Prop: generalised local homology} simplify the computation of $\text{H}_{\mathfrak{ a}}^{i}(M,N)$ (resp.  $\text{H}^{\mathfrak{ a}}_{i}(M,N)$) since the usual limits, $\underset{k}{\varinjlim}$ and $\underset{k}{\varprojlim}$ in the respective definitions get dropped when $N\in\mathfrak{R}^{M}_{\mathfrak{ a}}$ (resp. $N\in\mathfrak{ C}^{M}_{\mathfrak{ a}}$).
\end{paragraph}

\begin{paragraph}\noi Note that $\text{H}^{\mathfrak{ a}}_{i}(M,N)\ncong \text{H}^{\mathfrak{ a}}_{i}(N,M)$ in general; see for instance \cite[Remark 2.2 (iv)]{T T Nam-GLH for artinian modules}. Recall that an $R$-module $N$ is \textit{$\mathfrak{ a}$-adically complete} if $N\cong\Lambda_{ \mathfrak{a}}(N).$
	\end{paragraph}
\begin{prop}\label{Prop: Commutativity of GLH}
Let $M$ and $N$ be $R$-modules. If $M$ and $N$ are $\mathfrak{ a}$-adically complete $\mathfrak{ a}$-coreduced $R$-modules, then {\normalfont $\text{H}^{\mathfrak{ a}}_{i}(M,N)\cong \text{H}^{\mathfrak{ a}}_{i}(N,M).$}
\end{prop}
\begin{prf}
Since $M\cong \Lambda_{ \mathfrak{a}}(M)$ and $N\cong \Lambda_{ \mathfrak{a}}(N)$, it is enough to show that there exists an isomorphism
{\normalfont $\text{H}^{\mathfrak{ a}}_{i}(M,\Lambda_{ \mathfrak{a}}(N))\cong \text{H}^{\mathfrak{ a}}_{i}(N,\Lambda_{ \mathfrak{a}}(M))$} for every $\mathfrak{ a}$-coreduced $R$-modules $M$, $N$.
It follows from Proposition \ref{Lem: cor tensor with cor is cor} that $\Lambda_{ \mathfrak{a}}(N)\in\mathfrak{ C}^{M}_{\mathfrak{ a}}$, and $\Lambda_{ \mathfrak{a}}(M)\in\mathfrak{ C}^{N}_{\mathfrak{ a}}$. Moreover, $\Lambda_{ \mathfrak{a}}(M)\cong M/\mathfrak{ a}M$ and $\Lambda_{ \mathfrak{a}}(N)\cong N/\mathfrak{ a}N$ \cite[Proposition 2.3]{David-Application-I}. This data together with Proposition \ref{Prop: generalised local homology} and the symmetry property of the Tor functor, that is, $\text{Tor}^{R}_{i}(M,N)\cong \text{Tor}^{R}_{i}(N,M)$ for all $M,N\in R \text{-Mod}$ \cite[Theorem 7.1]{J.Rotman- Introduction hom.algebra} induces the isomorphisms 
$\text{H}^{\mathfrak{ a}}_{i}(M,\Lambda_{ \mathfrak{a}}(N))\cong \text{Tor}^{R}_{i}(M/\mathfrak{ a}M,\Lambda_{ \mathfrak{a}}(N))\cong \text{Tor}^{R}_{i}(M/\mathfrak{ a}M,N/\mathfrak{ a}N)$ $\cong \text{Tor}^{R}_{i}(N/\mathfrak{ a}N, M/\mathfrak{ a}M)\cong \text{H}^{\mathfrak{ a}}_{i}(N, \Lambda_{ \mathfrak{a}}(M)).$
\end{prf} 
\begin{prop}
Let $R$ be a ring, $S$ a multiplicatively closed subset of $R$, and $S^{-1}\mathfrak{ a}$ an ideal of the localized ring $S^{-1}R$. If $N\in\mathfrak{ C}_{\mathfrak{ a}}^{M}$, then for all $i\ge0$ 

 $$  S^{-1}\text{H}^{\mathfrak{ a}}_{i}(M,N)\cong \text{H}^{S^{-1}\mathfrak{ a}}_{i}(S^{-1}M, S^{-1}N).$$ 
 
\end{prop}

\begin{prf}	Suppose that $N\in\mathfrak{C}_{\mathfrak{ a}}^{M}$. By Proposition \ref{Prop: Localised modules}, $S^{-1}N\in \mathfrak{ C}^{S^{-1}M}_{S^{-1}\mathfrak{ a}}$ and  $\text{H}_{i}^{\mathfrak{ a}}(M,N)\cong \text{Tor}^{R}_{i}(M/\mathfrak{ a}M,N)$ [Proposition \ref{Prop: generalised local homology}]. This together with the fact that the Tor functor commutes with localization \cite[Proposition 7.17]{J.Rotman- Introduction hom.algebra}, \cite[Corollary 3.2.10]{weibel1995introduction} we get 
$$S^{-1}\text{H}^{\mathfrak{ a}}_{i}(M,N)\cong S^{-1}\text{Tor}^{R}_{i}(M/\mathfrak{ a}M,N)\cong \text{Tor}^{S^{-1}R}_{i}(S^{-1}(M/\mathfrak{ a}M), S^{-1}N)\cong $$ $$\text{Tor}^{S^{-1}R}_{i}(S^{-1}M/S^{-1} \mathfrak{ a}M, S^{-1}N)\cong \text{H}^{{S^{-1}\mathfrak{ a}}}_{i}(S^{-1}M, S^{-1}N)~\text{for all}~i\ge0.$$ 
\end{prf}
\begin{paragraph}\noi If $(R,\mathfrak{ m})$ is a local ring, and $M$ and $N$ are finitely generated $R$-modules, then $\text{H}^{\mathfrak{ m}}_{i}(M,N)$ is Artinian for all $i\ge 0$ \cite[Corollary 2.10]{T T Nam-GLH for artinian modules}.
\end{paragraph}
\begin{thm} \label{Prop: Artinianess + Noetherianess}
Let $\mathfrak{a}$ be an ideal of a ring $R$, and let {\normalfont $M,N\in R\text{-mod}$}.
\begin{enumerate}
	\item [$(1)$]  
  If $N\in \mathfrak{R}^{M}_{\mathfrak{ a}}$ (resp. $N\in\mathfrak{ C}^{M}_{\mathfrak{ a}}$), then {\normalfont $\text{H}_{\mathfrak{ a}}^{i}(M,N)$ (resp. $\text{H}^{\mathfrak{ a}}_{i}(M,N)$)} is Noetherian for all $i\ge0$.
  \item[$(2)$] Let  $R$ be Artinian. 
  If $N\in \mathfrak{R}^{M}_{\mathfrak{ a}}$ (resp. $N\in\mathfrak{ C}^{M}_{\mathfrak{ a}}$), then {\normalfont $\text{H}_{\mathfrak{ a}}^{i}(M,N)$ (resp. $\text{H}^{\mathfrak{ a}}_{i}(M,N)$)} is Artinian for all $i\ge0$.
\end{enumerate}
\end{thm}
\begin{prf}
 First note that $M/\mathfrak{ a}M$ is finitely generated, since $M$ is finitely generated. Suppose that $N\in\mathfrak{R}^{M}_{\mathfrak{ a}}$. Then $\text{H}^{i}_{\mathfrak{ a}}(M,N)
\cong \text{Ext}^{i}_{R}(M/\mathfrak{ a}M,N)$  [Proposition \ref{Prop: generalised local cohomology}], and {\normalfont $\text{Ext}^{i}_{R}(M/\mathfrak{ a}M,N)$} is finitely generated for all $i\ge0$ \cite[Theorem 7.36]{J.Rotman- Introduction hom.algebra}.
 It follows that {\normalfont $\text{H}^{i}_{\mathfrak{ a}}(M,N)$} is finitely generated for all $i$. Thus, {\normalfont $\text{H}^{i}_{\mathfrak{ a}}(M,N)$} is Noetherian for all $N\in\mathfrak{R}^{M}_{\mathfrak{ a}}$ and all $i\ge 0$. Moreover, by \cite[Theorem 7.20]{J.Rotman- Introduction hom.algebra}, $\text{Tor}^{R}_{i}(M/\mathfrak{ a}M,N)$ is finitely generated, and $\text{H}^{\mathfrak{ a}}_{i}(M,N)\cong \text{Tor}^{R}_{i}(M/\mathfrak{ a}M,N)$ for all $N\in\mathfrak{ C}^{M}_{\mathfrak{ a}}$ and all $i\ge0$ [Proposition \ref{Prop: generalised local homology}]. We conclude that $\text{H}^{\mathfrak{ a}}_{i}(M,N)$ over the Noetherian ring  $R$ is finitely generated for all $N\in\mathfrak{ C}^{M}_{\mathfrak{ a}}$ and all $i\ge0$ [Proposition \ref{Prop: generalised local homology}]. This proves part $(1)$. Similar proof applies to part $(2)$.
\end{prf}
\begin{cor}
	Let $\mathfrak{ a}$ be an ideal of a Noetherian (resp. an Artinian) ring $R$, and let {\normalfont $M,N\in R\text{-mod}$}.  If $N$ is $\mathfrak{ a}$-reduced (resp. $\mathfrak{ a}$-coreduced), then {\normalfont $\text{H}^{i}_{\mathfrak{ a}}(N)$} (resp. {\normalfont $\text{H}^{\mathfrak{ a}}_{i}(N)$)} is Noetherian (resp. Artinian) for all $i\ge0$.
	 
\end{cor}
\begin{prf}
This follows from Theorem \ref{Prop: Artinianess + Noetherianess}  with $M=R$.
\end{prf}  

\begin{prop}\label{Prop: coincidence of gen. local hom vs local cohomology}
	Let $\mathfrak{ a}$ be an ideal of  $R$, and let $N\in\mathfrak{ C}^{M}_{\mathfrak{ a}}$. Then  $~\text{for all}~ i\ge 0$,
{\normalfont$$\text{H}^{\mathfrak{ a}}_{i}(M,N)^{\vee}\cong \text{H}^{i}_{\mathfrak{ a}}(M,N^{\vee}).$$} 
\end{prop}
\begin{prf}
From the Ext-Tor Duality, we obtain an isomorphism $\text{Tor}^{R}_{i}(M/\mathfrak{ a}M, N)^{\vee}\cong \text{Ext}^{i}_{R}(M/\mathfrak{ a}M, N^{\vee})$ for all $i\ge0$. Since $N\in\mathfrak{ C}^{M}_{\mathfrak{ a}}$, and thus $N^{\vee}\in\mathfrak{R}^{M}_{\mathfrak{ a}}$ [Corollary \ref{Cor: X is cor iff its dual is red}], we conclude that $\text{H}^{\mathfrak{ a}}_{i}(M,N)^{\vee}\cong \text{H}^{i}_{\mathfrak{ a}}(M,N^{\vee})$ by Propositions \ref{Prop: generalised local cohomology} and  \ref{Prop: generalised local homology}.
\end{prf} 
\begin{paragraph}\noi Assume that $(R, \mathfrak{ m})$ is a local ring and {\normalfont $M\in R\text{-mod}$}. For all $i\ge 0$ and all {\normalfont $N\in R\text{-Mod}$},  $\text{H}^{\mathfrak{ a}}_{i}(M,N^{\vee})\cong \text{H}^{i}_{\mathfrak{ a}}(M,N)^{\vee}$ \cite[Proposition 2.3]{T T Nam-GLH for artinian modules}, \cite[Lemma 3.1]{T.T. Nam-generalised I-adic completion}.
\end{paragraph}
\begin{prop}\label{Prop: Duality_gen.(co)homology}
Let $\mathfrak{ a}$ be an ideal of a ring  $R$, and {\normalfont $M,N\in R\text{-Mod}$} with $M$ finitely generated. 
For each $N\in\mathfrak{R}^{M}_{\mathfrak{ a}}$ and each $i\ge 0$, {\normalfont $$ \text{H}^{\mathfrak{ a}}_{i}(M,N^{\vee})\cong \text{H}^{i}_{\mathfrak{ a}}(M,N)^{\vee}.$$} 
\end{prop}
\begin{prf}
By \cite[Definition 3.3.9]{weibel1995introduction}, $M$ has a resolution $\textbf{F}\rightarrow M$ such that each $\text{F}^{i}\in\textbf{F}$ is a finitely generated free $R$-module. It follows that  $\text{Tor}^{R}_{i}(M,N^{\vee})\cong \text{Ext}^{i}_{R}(M,N)^{\vee}$ for all $i$ \cite[1.4.8]{Peter-Schenzez}. Since $N\in \mathfrak{R}^{M}_{\mathfrak{ a}}$, and thus $N^{\vee}\in\mathfrak{ C}^{M}_{\mathfrak{ a}}$ [Corollary \ref{Cor: X is cor iff its dual is red}],  $\text{H}^{\mathfrak{ a}}_{i}(M,N^{\vee})\cong \text{Tor}^{R}_{i}(M/\mathfrak{a}M,N^{\vee})$, see Proposition \ref{Prop: generalised local homology}. These data, and Proposition \ref{Prop: generalised local cohomology} yield
$\text{H}^{\mathfrak{ a}}_{i}(M,N^{\vee})\cong \text{Tor}^{R}_{i}(M/\mathfrak{ a}M, N^{\vee})\cong \text{Ext}^{i}_{R}(M/\mathfrak{ a}M,N)^{\vee}\cong \text{H}^{i}_{\mathfrak{ a}}(M,N)^{\vee}$ for all $i\ge 0$.
\end{prf}  
\section{Application via spectral sequences}
\begin{paragraph}\noi
	In this section some  Grothendieck's (co)homological spectral sequences and their applications to the generalised local (co)homology modules in the setting  of $\mathfrak{ a}$-reduced modules and $\mathfrak{ a}$-coreduced modules with respect to $M$ are established. In particular,  $\mathfrak{R}^{M}_{\mathfrak{ a}}$ and $\mathfrak{ C}^{M}_{\mathfrak{ a}}$ respectively (1) facilitate the vanishing of $\text{H}^{p}_{\mathfrak{ a}}(M,\text{H}^{q}_{\mathfrak{ a}}(M,N))$ and $\text{H}^{\mathfrak{ a}}_{p}(M,\text{H}^{\mathfrak{ a}}_{q}(M,N))$ whenever either $p\ne0$ or $q\ne0$; and (2) determine when $\text{H}^{q}_{\mathfrak{ a}}(M,N)$ (resp. $\text{H}_{q}^{\mathfrak{ a}}(M,N)$) belong to $\mathfrak{R}^{M}_{\mathfrak{ a}}$ (resp. $\mathfrak{C}^{M}_{\mathfrak{ a}}$). 
As in section six, all the rings in this section are Noetherian.
\end{paragraph}
 
\begin{prop}\label{Prop: spectral sequences} 
	Let  $R$ be a ring, {\normalfont $M,N\in R \text{-Mod}$}, and let $p,q\in\Z$.  There exist Grothendieck's (co)homological spectral sequences: 
	\begin{enumerate}
		\item [$(1)$] {\normalfont $\text{E}_{2}^{p,q}=\text{H}^{p}_{\mathfrak{a}}(\text{Ext}^{q}_{R}(M,N))\Rightarrow \text{H}^{p+q}_{\mathfrak{a}}(M,N).$} 
		\item[$(2)$] 
		$ {\normalfont \text{E}_{2}^{p,q}= \text{Ext}^{p}_{R}(M,\text{H}^{q}_{\mathfrak{a}}(N))\Rightarrow \text{H}^{p+q}_{\mathfrak{a}}(M,N)} ~~ \text{provided}~{\normalfont M \in R \text{-mod}}. $
		\item[$(3)$] ${\normalfont \text{E}^{2}_{p,q}=  \text{H}^{\mathfrak{a}}_{p}(\text{Tor}^{R}_{q}(M,N))\Rightarrow \text{H}^{\mathfrak{a}}_{p+q}(M,N)} ~\text{for all}~{\normalfont M, N\in R \text{-mod}}.$
		\item[$(4)$] $ {\normalfont \text{E}^{2}_{p,q}= \text{Tor}_{p}^{R}(M,\text{H}_{q}^{\mathfrak{a}}(N))\Rightarrow \text{H}_{p+q}^{\mathfrak{a}}(M,N)} ~~ \text{for all}~{\normalfont M, N\in R \text{-mod}}.$
	\end{enumerate}
\end{prop}
\begin{prf} The first two parts are in \cite[Proposition 2.1]{Freitas-Schenzel GLD}. We prove part $(3)$, since the proof of part $(4)$ is similar.
  
Let $P$ be a projective $R$-module and $i\ne0$. We have
$\text{L}_{i}\Lambda_{\mathfrak{a}}(M\otimes_{R}P):=\text{H}^{\mathfrak{a}}_{i}(M,P):=\underset{t}\varprojlim\text{Tor}^{R}_{i}(M/\mathfrak{a}^{t}M,P)=0.$ So, by \cite[Corollary 5.8.4 ]{weibel1995introduction},	$\text{E}^{2}_{p,q}=\text{L}_{p}\Lambda_{\mathfrak{a}}(\text{L}_{q}(M\otimes_{R}N))\Rightarrow \text{L}_{p+q}\Lambda_{\mathfrak{a}}(M\otimes_{R}N)~~\text{for all} ~N\in R\text{-mod}$. This shows that
		$\text{H}^{\mathfrak{a}}_{p}(\text{Tor}^{R}_{q}(M,N))\Rightarrow \text{H}_{p+q}^{\mathfrak{a}}(M,N) ~~\text{for all}~ N\in R\text{-mod}~\text{and}~\text{all}~ p,q\in\Z.$
 
\end{prf}
 
\begin{paragraph}\noi
 For any $M\in R\text{-Mod}$, let $\mathfrak{B}_{\mathfrak{ a}}^{M}:= \mathfrak{R}_{\mathfrak{ a}}^{M}\cap \mathfrak{C}_{\mathfrak{ a}}^{M}$.  
	\end{paragraph}
\begin{cor}\label{ Cor--Noeth. examples in GLH,GLC}
 Let  $R$ be a ring, and {\normalfont $M,N\in R\text{-mod}$}. For all $p\in\Z$ and all  $\mathfrak{a}$-coreduced $M$, {\normalfont $\text{H}^{p}_{\mathfrak{a}}(M,N), \text{H}^{\mathfrak{ a}}_{p}(M,N)\in \mathfrak{B}^{M}_{\mathfrak{ a}}$}. 
\end{cor}
\begin{prf}
   $\text{H}^{p}_{\mathfrak{a}}(\text{Hom}_{R}(M,N))\cong \text{H}^{p}_{\mathfrak{a}}(M,N)$ and $\text{H}^{\mathfrak{ a}}_{p}(M\otimes_{R}N)\cong \text{H}^{\mathfrak{ a}}_{p}(M,N)$ for all $p\in\Z$ and all $N\in R\text{-mod}$ by Proposition \ref{Prop: spectral sequences} part $(1)$ (resp. part $(3)$).  Since $M$ is $\mathfrak{ a}$-coreduced, $\text{Hom}_{R}(M,N)$ is $\mathfrak{ a}$-reduced \cite[Proposition 2.6 (1)]{David-Application-I} and $M\otimes_{R}N$ is $\mathfrak{ a}$-coreduced [Proposition \ref{Lem: cor tensor with cor is cor}]. By \cite[Corollary 6.1]{ssevviiri: App II}, $\text{H}^{p}_{\mathfrak{ a}}(\text{Hom}_{R}(M,N), \text{H}^{\mathfrak{ a}}_{p}(M\otimes_{R}N)\in\mathfrak{B}_{\mathfrak{ a}}^{M}$.
   Therefore, $\text{H}^{p}_{\mathfrak{ a}}(M,N), \text{H}^{\mathfrak{ a}}_{p}(M,N) \in\mathfrak{B}_{\mathfrak{ a}}^{M}$. 
\end{prf} 
\begin{prop}\label{Prop: Spectral sequences}
	Let $\mathfrak{a}$ be an ideal of $R$, and let {\normalfont $M,N\in R\text{-Mod}$}. 
	\begin{enumerate}
		\item [$(1)$] If $N\in\mathfrak{R}^{M}_{\mathfrak{ a}}$, then there is a Grothendieck's cohomological spectral sequence {\normalfont $\text{E}^{p,q}_{2}=\text{H}_{\mathfrak{a}}^{p}(\text{Ext}^{q}_{R}(M,N))\Rightarrow  \text{Ext}^{p+q}_{R}(M/\mathfrak{a}M, N)~~\text{for all}~ p,q\in\Z.$}
		\item[$(2)$]Let {\normalfont $M,N\in R\text{-mod}$}.  If $N\in\mathfrak{ C}^{M}_{\mathfrak{ a}}$, then there exists a homological spectral sequence
		{\normalfont 	$\text{E}_{p,q}^{2}=\text{H}^{\mathfrak{a}}_{p}(\text{Tor}^{q}_{R}(M,N))\Rightarrow  \text{Tor}_{p+q}^{R}(M/\mathfrak{a}M, N)~~\text{for all}~ p,q\in\Z.$}
	\end{enumerate}
\end{prop}
\begin{prf}  Since by Proposition \ref{Prop: generalised local cohomology}, $\text{H}^{i}_{\mathfrak{ a}}(M,N)\cong \text{Ext}^{i}_{R}(M/\mathfrak{ a}M,N)$ for all $i\ge0$, the conclusion for part $(1)$ follows by Proposition \ref{Prop: spectral sequences} $(1)$. For part $(2)$, see Proposition \ref{Prop: generalised local homology} and Proposition \ref{Prop: spectral sequences} $(3)$.
\end{prf}
\begin{prop}\label{Prop: Reduced and Coreduced vs GLH}
	Let $\mathfrak{a}$ be an ideal of a $R$, and let {\normalfont $M,N\in R\text{-Mod}$} with $M$ finitely generated. Let $q\in\Z$. 
	\begin{enumerate}
		\item [$(1)$] If {\normalfont $\text{H}^{q}_{\mathfrak{a}}(N)\in\mathfrak{R}^{M}_{\mathfrak{ a}}$}, then {\normalfont $\text{H}^{q}_{\mathfrak{a}}(M,N)\in \mathfrak{R}^{M}_{\mathfrak{ a}}$}.
		\item[$(2)$] If {\normalfont $\text{H}_{q}^{\mathfrak{a}}(N)\in \mathfrak{C}^{M}_{\mathfrak{ a}}$}, then {\normalfont $\text{H}_{q}^{\mathfrak{a}}(M,N)\in \mathfrak{C}^{M}_{\mathfrak{ a}}$} for all {\normalfont $N\in R\text{-mod}$}.
	\end{enumerate}
\end{prop}
\begin{prf}
			By Proposition \ref{Prop: spectral sequences} $(2)$, $\text{Hom}_{R}(M,\text{H}^{q}_{\mathfrak{a}}(N))\cong \text{H}^{q}_{\mathfrak{a}}(M,N)$ follows for all $q$ and all $N\in R\text{-Mod}$. Since $\text{H}^{q}_{\mathfrak{a}}(N)\in\mathfrak{R}^{M}_{\mathfrak{ a}}$, we have $\text{Hom}_{R}(M,\text{H}^{q}_{\mathfrak{a}}(N))\in\mathfrak{R}^{M}_{\mathfrak{ a}}$ [Proposition \ref{Prop: Lemma 2.2 of our first paper}]. So, $\text{H}^{q}_{\mathfrak{a}}(M,N) \in\mathfrak{R}^{M}_{\mathfrak{ a}}$. This establishes part $(1)$. Part $(2)$ follows from Proposition \ref{Prop: Lemma 2.2 of our first paper} and Proposition \ref{Prop: spectral sequences} $(4)$, similarly.
 
\end{prf}

\begin{paragraph}\noi We conclude this section by presenting Proposition \ref{Gen.Prop 6.1 _D2} and Proposition \ref{Gen.Prop 6.2 _D2} as the generalised versions of \cite[Proposition 6.1]{ssevviiri: App II} and \cite[Proposition 6.2]{ssevviiri: App II} respectively. 
	\end{paragraph}
 
\begin{prop}\label{Gen.Prop 6.1 _D2}
	Let $R$ be a von-Neumann regular ring, and let $\mathfrak{ a}$ be an ideal of $R$. For all $p,q\in\Z$ and all {\normalfont $M,N\in R\text{-mod}$},
	
    $$ {\normalfont \text{H}^{\mathfrak{ a}}_{p}(M,\text{H}^{\mathfrak{ a}}_{q}(M,N))}= \begin{cases} 
		\Lambda_{ \mathfrak{a}}(M, \Lambda_{ \mathfrak{a}}(M,N), & p=q=0,\\
		0 & \text{otherwise}.	
	\end{cases}$$ and the associated Grothendieck's spectral sequence is given by
{\normalfont $$\text{E}^{2}_{0,0}= \Lambda_{ \mathfrak{a}}(M,\Lambda_{ \mathfrak{a}}(M,N))\Rightarrow M/\mathfrak{ a}M\otimes_{R}M/\mathfrak{ a}M\otimes_{R}N;$$} and
{\normalfont $$\text{E}^{2}_{p,q}= {\normalfont \text{H}^{\mathfrak{ a}}_{p}(M,\text{H}^{\mathfrak{ a}}_{q}(M,N))}\Rightarrow 0 ~\text{if either}~ p\ne0 ~\text{or}~ q\ne0.$$}
\end{prop}
\begin{prf}
 For a von-Neumann regular ring $R$, every $R$-module is $\mathfrak{ a}$-coreduced with respect to $M$. So, by Proposition \ref{Prop: generalised local homology}, 
  $$ {\normalfont \text{H}^{\mathfrak{ a}}_{p}(M,\text{H}^{\mathfrak{ a}}_{q}(M,N))}\cong 
  	\text{Tor}_{p}^{R}(M/\mathfrak{ a}M, \text{Tor}^{R}_{q}(M/\mathfrak{ a}M, N))$$
$$= \begin{cases} 
	M/\mathfrak{ a}M \otimes_{R} M/\mathfrak{ a}M\otimes_{R} N, & p=q=0\\
	0 & \text{otherwise}.	
\end{cases}$$
This is because if either $p\ne0$ or $q\ne0$, then $\text{Tor}_{p}^{R}(M/\mathfrak{ a}M, \text{Tor}^{R}_{q}(M/\mathfrak{ a}M, N))=0$ as modules over von-Neumann regular rings are flat. By Proposition \ref{sec1:Prop- strongly I-coreduced}, $M/\mathfrak{ a}M\otimes_{R}M/\mathfrak{ a}M\otimes_{R}N\cong M/\mathfrak{ a}M\otimes_{R}\Lambda_{ \mathfrak{a}}(M,N)\cong\Lambda_{ \mathfrak{a}}(M,\Lambda_{\mathfrak{ a}}(M,N))$. So, 
$$ {\normalfont \text{H}^{\mathfrak{ a}}_{p}(M,\text{H}^{\mathfrak{ a}}_{q}(M,N))}\cong \begin{cases} 
	\Lambda_{ \mathfrak{a}}(M, \Lambda_{ \mathfrak{a}}(M, N)), & p=q=0\\
	0 & \text{otherwise}.	
\end{cases}$$ 

\end{prf}

\begin{prop}\label{Gen.Prop 6.2 _D2}
	Let $R$ be an Artinian  von-Neumann regular ring, and let $\mathfrak{ a}$ be an ideal of $R$. If $M$ is a finitely generated $R$-module, then 
	
	$$ {\normalfont \text{H}_{\mathfrak{ a}}^{p}(M,\text{H}_{\mathfrak{ a}}^{q}(M,N))}= \begin{cases} {\normalfont
		\Gamma_{ \mathfrak{a}}(M, \Gamma_{ \mathfrak{a}}(M,N))} & p=q=0,\\
		0 & \text{otherwise}	
	\end{cases}$$ 
for all $p,q\in\Z$ and all {\normalfont $N\in R\text{-Mod}$}. The associated Grothendieck's spectral sequence is given by
{\normalfont $$E^{2}_{0,0}= {\normalfont \text{H}_{\mathfrak{ a}}^{0}(M,\text{H}_{\mathfrak{ a}}^{0}(M,N))}\cong \Gamma_{ \mathfrak{a}}(M,\Gamma_{ \mathfrak{a}}(M,N))\Rightarrow \text{Hom}_{R}(M/\mathfrak{ a}M, \text{Hom}_{R}( M/\mathfrak{a}M, N),$$} and
$$\text{E}^{2}_{p,q}= {\normalfont \text{H}_{\mathfrak{ a}}^{p}(M,\text{H}_{\mathfrak{ a}}^{q}(M,N))}\Rightarrow 0, \text{if either}~ p\ne0 ~\text{or}~ q\ne0.$$
\end{prop}
\begin{prf} 
 Every $R$-module over a von-Neumann regular ring is $\mathfrak{ a}$-reduced with respect to $M$. So, by Proposition \ref{Prop: generalised local cohomology}, 
		$$ {\normalfont \text{H}_{\mathfrak{ a}}^{p}(M,\text{H}_{\mathfrak{ a}}^{q}(M,N))}\cong 
			\text{Ext}^{p}_{R}(M/\mathfrak{ a}M,\text{Ext}_{R}^{q}(M/\mathfrak{ a}M, N)).$$ 
	Since $R$ is Artinian von-Neumann regular, $R$ is semi-simple and therefore every $R$-module is projective. Thus,  $\text{Ext}^{p}_{R}(M/\mathfrak{ a}M, \text{Ext}^{q}_{R}(M/\mathfrak{ a}M,N))=0$ provided either $p\ne 0$ or $q\ne 0$. For $p=q=0$, $$\text{Ext}^{p}_{R}(M/\mathfrak{ a}M,\text{Ext}^{q}_{R}(M/\mathfrak{ a}M, N)= \text{Hom}_{R}(M/\mathfrak{ a}M,\text{Hom}_{R}( M/\mathfrak{a}M, N)).$$ Furthermore, since every $R$-module over a von-Neumann regular ring is $\mathfrak{ a}$-reduced with respect to $M$, $\text{Hom}_{R}(M/\mathfrak{ a}M,\text{Hom}_{R}( M/\mathfrak{a}M, N))\cong\text{Hom}_{R}(M/ \mathfrak{ a}M, \Gamma_{ \mathfrak{a}}(M,N))\cong\Gamma_{ \mathfrak{a}}(M,\Gamma_{\mathfrak{ a}}(M,N))$ by Proposition \ref{Prop: reduced w.r.t M}. Therefore, 
		$$ {\normalfont \text{H}_{\mathfrak{ a}}^{p}(M,\text{H}_{\mathfrak{ a}}^{q}(M,N))}\cong \begin{cases} 
			\Gamma_{ \mathfrak{a}}(M, \Gamma_{ \mathfrak{a}}(M, N)), & p=q=0\\
			0 & \text{otherwise}.	
		\end{cases}$$
	
	\end{prf}
 
 \subsection*{Acknowledgement} 
 \begin{paragraph}\noi
 	We acknowledge support from the Eastern Africa Algebra Research Group (EAALG), the EMS-Simons for Africa Program and the International Science Program (ISP).  
 \end{paragraph}

\end{document}